\documentclass[11pt]{article}
\usepackage{amsmath}\usepackage{amssymb}\usepackage{amscd}
\usepackage{color}\usepackage{tikz} 
\usepackage{stackrel}

\newenvironment{annotacia}{\centerline{\sc Abstract}\vspace{2mm}\narrower\narrower\sf}

\def\theequation{\thesection.\arabic{equation}}
\makeatletter\@addtoreset{equation}{section}\makeatother

\def\nn{\nonumber}\def\lb{\label}
\def\Mt{M \raisebox{1mm}{$\intercal$}}
\def\be{\begin{equation}}
\def\ee{\end{equation}}
\def\ba{\begin{eqnarray}}
\def\ea{\end{eqnarray}}
\def\tr{{\rm Tr}\,}
\def\Tr#1{{\rm Tr}_{\! R^{\mbox{\scriptsize$(#1)$}}}}
\def\TR#1#2{{\rm Tr}_{\! #2^{\mbox{\,\scriptsize$(#1)$}}}}
\def\str#1{\rule[#1mm]{0pt}{1mm}}

\newcounter{theorem}\makeatletter
\@addtoreset{theorem}{section}\makeatother

\newtheorem{prop}[theorem]{Proposition}
\newtheorem{rem}[theorem]{Remark}
\newtheorem{lem}[theorem]{Lemma}
\newtheorem{def-lem}[theorem]{Definition-Lemma}
\newtheorem{def-prop}[theorem]{Definition-Proposition}
\newtheorem{defin}[theorem]{Definition}
\newtheorem{theor}[theorem]{Theorem}

\begin{document}
\begin{center}
{\Large \textbf{Cayley-Hamilton Theorem for Symplectic\\[4pt] Quantum Matrix Algebras}}

\vspace{1cm} {\large \textbf{Oleg Ogievetsky$^{\diamond\,\dag\,\ddagger}$ and Pavel Pyatov$^{\ast\,\star
}$}}

\vskip .8cm $^{\diamond}$Aix Marseille Universit\'{e}, Universit\'{e} de
Toulon, CNRS, \\ CPT UMR 7332, 13288, Marseille, France

\vskip .3cm $^{\dag}${I.E.Tamm Department of Theoretical Physics, P.N. Lebedev Physical Institute, Leninsky prospekt 53, 119991 Moscow, Russia}

\vskip .3cm $^{\ddagger}${Inst. of the Information Transmission Problems, RAS,
Moscow, Russia}

\vskip .3cm $^{\ast}${HSE University, 20 Myasnitskaya street, Moscow
101000, Russia}

\vskip .3cm $^{\star}${Bogoliubov Laboratory of Theoretical Physics, JINR, 141980 Dubna, Moscow region, Russia}
\end{center}

\vskip 1cm 
\begin{annotacia} \noindent
We  establish the analogue of the Cayley--Hamilton theorem for the quantum matrix algebras of the symplectic
type.
 \end{annotacia}

\newpage

\tableofcontents

\bigskip\bigskip\bigskip

\section{Introduction}\lb{sec1}

\smallskip
Let $V$ be a vector space equipped with a bilinear nondegenerate (symmetric or antisymmetric) form.
The Brauer algebra \cite{Br} generalizes the tower of the centralizing algebras which appears in the
 Brauer--Schur--Weyl duality, related to $V$. 
The Birman--Murakami--Wenzl algebra \cite{BW,M1} is the quantum deformation of the Brauer algebra. Important 
particular cases of local representations of the tower of the 
Birman--Murakami--Wenzl (BMW) algebra are constructed with the use of orthogonal and 
symplectic R-matrices. These R-matrices give rise to the quantum matrix algebras of the
orthogonal and symplectic types. More precisely, a quantum matrix algebra is defined by a compatible pair 
$\{ R, F\}$ of  R-matrices, we recall the definitions below. The general structure properties of quantum matrix algebras 
with $R$ of the BMW type were investigated in \cite{OP}. 
In the present work we mainly assume $R$ to be 
of the symplectic type. Our principal goal is to derive, for the quantum matrices of the symplectic type,
an analogue of the Cayley--Hamilton identity and to use it for a
description of the spectra of the corresponding quantum matrices.

\smallskip
In Section \ref{BMWaRm} we recall the necessary facts about the Birman-Murakami-Wenzl (BMW)
algebras, their R-matrix realizations, specializations to the symplectic and orthogonal cases and some 
R-matrix technique. 

\vskip .1cm
Section \ref{secqma} contains the information from \cite{OP} about the quantum matrix algebra, 
its `characteristic subalgebra' (the subalgebra to which the coefficients of the Cayley-Hamilton identity belong)
and the $\star$-multiplication.

The main results are  in Section \ref{secCHSp}. Here we establish the Cayley--Hamilton theorem for the 
symplectic quantum matrix algebras. Classically, the symplectic group is defined by the
condition $M^t \Omega M = \Omega$ where $ \Omega$ is the symplectic form. However we have to a work with a 
bigger group defined by the condition $M^t \Omega M = g\Omega$ where $g$ is a constant, that is with the
group of transformations which preserve the form up to a
multiplicative factor.  We  call `2-contraction' the quantum analogue of this factor. It is an element $g$ of 
the quantum matrix algebra.
For a general compatible pair $\{ R, F\}$, the element  $g$ is not necessarily central so
we cannot harmlessly set it to 1. We establish a strengthened form of the Cayley--Hamilton theorem which does 
not assume the invertibility of the element $g$ (and which is equivalent to the Cayley--Hamilton theorem under the
assumption of the invertibility of $g$). 

\vskip .1cm
Next, we define in Section \ref{secCHSp} a homomorphism from the
characteristic subalgebra to the algebra of symmetric polynomials in some set of commuting
(``spectral") variables. 
The nature of this homomorphism reflects the reciprocity properties of the characteristic
polynomials for the symplectic matrices. 
The Cayley-Hamilton
identities under the action of this homomorphism are completely factorized  and hence the
spectral variables can be treated as eigenvalues of the quantum matrix. 
We then give the spectral parameterization of the  three series of 
elements of the characteristic subalgebra: the power sums $p_i$, the elementary symmetric functions $a_i$ and 
the complete symmetric functions $s_i$.

\vskip .1cm
Section \ref{secCHSp} contains also the low-dimensional examples illustrating the Cayley--Hamilton theorem
for two most known quantum matrix algebras: the algebra of functions on the quantum group (corresponding to 
the compatible pair $\{ R, P\}$ where $P$ is the flip) and the reflection equation algebra \cite{C,KS} 
(the reflection equation algebra corresponds to 
the compatible pair $\{ R, R\}$). Also we discuss the classical limit of the Cayley--Hamilton theorem.

\vskip .1cm
The Cayley--Hamilton identity for the quantum matrix algebras of the orthogonal type will be considered in a separate publication.

\section{BMW algebra and their R-matrix representations}\lb{BMWaRm}

In this section we present definitions and describe necessary facts about the Birman-Murakami-Wenzl algebras and the BMW type R-matrices. We follow notation of ref. \cite{OP} where the reader can find detailed
derivations and the references. Later in the section we investigate two families of the BMW type R-matrices, the $Sp(2k)$ type and the $O(k)$ type R-matrices. They are related, respectively, to the symplectic and orthogonal series of the quantum groups. We identify particular conditions on the eigenvalues which are specific for these families of R-matrices. In the following section we will use the symplectic R-matrices for the definition of  $Sp(2k)$ type quantum matrix algebras. Specific properties of the
$Sp(2k)$ type R-matrices will then dictate a form of the Cayley-Hamilton identites in these algebras.

\subsection{BMW algebra}

The {\em Birman-Murakami-Wenzl (BMW) algebra} ${\cal W}_{n}(q,\mu)$ \cite{BW,M1} depending an two complex parameters $q\in{\mathbb{C}}\backslash \{0,\pm 1\}$ and $\mu\in{\Bbb C}\backslash\{0,q,-q^{-1}\}$
is defined in terms of generators $\{\sigma_i, \kappa_i\}_{i=1}^{n-1}$ and  relations
\ba
\nn
\sigma_i \sigma_{i+1} \sigma_i \, =\,  \sigma_{i+1} \sigma_i \sigma_{i+1},&&
\sigma_i \sigma_j \, =\,  \sigma_j \sigma_i\qquad  \qquad\;\; \forall\; i,j:\; |i-j|>1 ,
\\[2pt]
\lb{kappa}
\sigma_i \kappa_i \, =\,  \kappa_i \sigma_i \,=\, \mu \kappa_i ,
&&
\kappa_{i} = { \textstyle {(q1-\sigma_i)(q^{-1} 1 +\sigma_i)\over \mu (q-q^{-1})}} , 
\\[2pt]
\nn
\kappa_{i+1}\kappa_i \, =\, \kappa_{i+1}\sigma^{\pm 1}_{i}\sigma^{\pm 1}_{i+1}, 
&&
\kappa_i\kappa_{i+1}\kappa_i \, =\,\kappa_i\;\;\;\;\;\;\qquad \forall\; i .
\ea
The first line is the Artin's presentation of the braid group ${\cal B}_n$; the rest of
relations define the quotient algebra ${\cal W}_{n}(q,\mu)\subset {\Bbb C}[{\cal B}_{n}]$.
\smallskip

Imposing further restrictions on the parameters
\be
\lb{mu}
j_q :={ q^j - q^{-j}\over q-q^{-1}}\neq\, 0\, , \quad \mu \neq \mp\, q^{\mp(2j-3)} \quad \forall\;j = 2,3,\dots , n\ ,
\ee
one can define recursively two sets of idempotents
\ba
\lb{ind1}
a^{(1)} := 1,\hspace{3.9cm} &&\;\;\, s^{(1)}\, :=\, 1,
\\[2pt]
\lb{a^k}
a^{(i+1)}\, :=\,{q^i\over (i+1)_q} a^{(i)}\, \sigma^{-}_i(q^{-2i})\,
a^{(i)} , &&
s^{(i+1)} :={q^{-i}\over (i+1)_q}\, s^{(i)}\, \sigma^{+}_i(q^{2i})\,
s^{(i)},
\ea
where
$$
\sigma_i^{\pm}(x)\, :=\,  1\, +\, {x-1\over q-q^{-1}}\, \sigma_i\, +\,
{\mu(x-1)\over \mu\mp q^{\mp 1}  x}\, \kappa_i\, .
$$
The idempotents $a^{(n)}$ and $s^{(n)}$ in the algebra
${\cal W}_{n}(q,\mu)$ are primitive. They correspond  to the $q$-deformations of the `trivial' ($\sigma_i\mapsto q$) and the `alternating' ($\sigma_i\mapsto -q^{-1}$)
one-dimensional representations. 
Therefore, they are called an {\em $n$-th order antisymmetrizer} and
an {\em $n$-th order symmetrizer}, respectively.

\subsection{R-matrices and their compatible pairs}

Let $V$ denote a finite dimensional ${\Bbb C}$-linear space, $\dim V = \mbox{\sc n}$. Fixing some
basis $\{v_i\}_{i=1}^{\mbox{\footnotesize \sc n}}$ in $V$ we identify elements $X\in
{\rm End}(V^{\otimes n})$ with  matrices $X_{i_1 i_2 \dots i_n}^{j_1 j_2 \dots j_n}$.

\vskip .1cm
Let $X\in {\rm End}(V^{\otimes k})$, $k\leq n$. For $1\leq m \leq n-k+1$, denote by $X_m\in {\rm End}(V^{\otimes n})$ an operator 
given by the matrix
$$
(X_m)_{i_1 \dots i_n}^{j_1 \dots j_n}\ :=\ I_{i_1\dots i_{m-1}}^{j_1\dots j_{m-1}}\
X_{i_m\dots i_{m+k-1}}^{j_m\dots j_{m+k-1}}\ I_{i_{m+k}\dots i_n}^{j_{m+k}\dots j_n}\ .
$$
Here $I$ denotes the identity operator.

\smallskip
An element $R\in {\rm Aut}(V^{\otimes 2})$ that fulfills 
an equation
$$
R_{1}\, R_{2}\, R_{1}\, = \, R_{2}\, R_{1}\, R_{2}\ .
$$
is called an {\em R-matrix}.
The permutation operator $P$, defined by $P(u\otimes v)= v\otimes u \;\;\; \forall\; u, v\in V\,$ is the R-matrix. The
operator $R^{-1}$ is the R-matrix iff $R$ is. 

\vskip .1cm
Any R-matrix $R$ generates
representations $\rho_R$ of the series of braid groups ${\cal B}_n$, $n=2,3,\dots$
$$
\rho_R:\, {\cal B}_n\rightarrow {\rm Aut}(V^{\otimes n})\ ,\quad
\sigma_i \mapsto   R_i, \quad 1\leq i\leq n-1 .
$$
An R-matrix is called {\it skew invertible} if there exists an operator
${\Psi_{_{\hspace{-1.mm}R}}}\in {\rm End}(V^{\otimes 2})$ such that
\be
\tr_{\!(2)} R_{12} {\Psi_R}_{23} =\tr_{\!(2)} {\Psi_R}_{12} R_{23} = P_{13}\, .
\lb{s-inv}
\ee
Here we use notation $X_{ij}$ which shows explicitly indices $i$ and $j$ of the spaces where
operator $X$ acts, e.g., $P_{13}= P_{i_1 i_3}^{j_1 j_3}\, I_{i_2}^{j_2}$. Symbol $\tr_{\!(i)}$ means taking the trace in the vector space with index $i$.

\vskip .2cm
For a skew invertible R-matrix $R$ define the operator $D_R\in \mbox{End}(V)$
\be
\lb{CandD}
(D_R)_1:={\rm Tr}_{(2)}{\Psi_R}_{12} .
\ee
The operator $R$ is called {\em strict skew invertible} if  $D_R$ is invertible. The R-matrix $R^{-1}$ is skew invertible iff $R$ is strict skew invertible \cite{I,O}, the corresponding operator $D_{R^{-1}}$ reads
$$
(D_{R^{-1}})_2= \left( {\rm Tr}_{(1)}{\Psi_R}_{12}\right)^{-1}.
$$

With a skew invertible R-matrix $R$ we associate a linear map on the space of $\mbox{\sc n}\!\times\! \mbox{\sc n}$ matrices whose entries belong to some $\mathbb{C}$-linear space $W$
$$
{\rm Tr\str{-1.3}}_R:\; {\rm End}(V)\otimes W\,\rightarrow \, W ,\ 
{\rm Tr\str{-1.3}}_R(M) ={\textstyle \sum_{i,j=1}^{\mbox{\footnotesize\sc n}}}{(D_R)}_i^jM_j^i\, , \ 
M\in{\rm End}(V)\otimes W\,,
$$
This map is called an  {\em R-trace}. 

\smallskip
It is easy to check that the R-matrix $P$ is strict skew invertible and $\tr_{\! P}$ coincides with the usual trace.  A characteristic property of the R-trace map is
\be
\lb{Rtr}
\Tr{2} R_1 = I_1.
\ee

An ordered pair $\{ R, F\}$ of two R-matrices $R$ and $F$  is
called {\em a compatible R-matrix pair} if the following conditions
\be 
R_1\, F_2\, F_1\, =\, F_2\, F_1\, R_2\, ,\qquad R_2\, F_1\, F_2\, =\, F_1\, F_2\, R_1\, ,
\lb{sovm}
\ee
are satisfied. The equalities (\ref{sovm}) are called {\em twist relations}. Clearly, $\{ R,P\}$ and $\{ R,R\}$ are compatible pairs of R-matrices.

\vskip .1cm
A compatible pair of R-matrices $\{R,F\}$ gives rise to a new R-matrix
\be
\lb{R_f}
R_f := F^{-1} R F\, ,
\ee
called the {\em twisted} R-matrix. The R-matrix pair $\{ R_f , F\}$ is compatible. If $R$ is skew invertible and $F$ is strict skew invertible, then $R_f$ is skew invertible; if additionally $R$ is strict skew invertible, then $R_f$ is strict skew invertible as well \cite{OP}.

\subsection{BMW type R-matrices}

Assume that an R-matrix $R$ satisfies a third order minimal
characteristic polynomial 
\be
\lb{charR}
(qI-R)(q^{-1}I+R)(\mu I-R)=0 ,
\ee
and an element
\be
\lb{K}
K := \mu^{-1} (q-q^{-1})^{-1}\, (qI-R)(q^{-1}I+R)\ee
fulfills  conditions
\be
\lb{bmwRa}
K_2\, K_1 \, = \,   K_2\, R_1^{\pm 1}\, R_2^{\pm 1}, 
\qquad
K_1\, K_2\, K_1 \, = \, K_1.
\ee
In this case  $R$ generates representations 
$\rho_R$ of the tower of the BMW algebras ${\cal W}_n(q,\mu)
\rightarrow {\rm End}(V^{\otimes n})$ $\forall\, n>1$
$$
\rho_R:\,{\cal W}_n(q,\mu)\rightarrow {\rm Aut}(V^{\otimes n})\ ,\quad
\sigma_i \mapsto   R_i,\quad \kappa_i \mapsto K_i,\quad 1\leq i\leq n-1 .
$$
Such R-matrix is said to be of {\em BMW type}.
\smallskip

If the R-matrix $R$ is skew invertible and of the BMW type, then it is strict skew invertible and the rank of the associated operator $K$ (\ref{K}) equals 1; the R-trace map in this case fulfills equalities \cite{IOP3}
\be
\lb{BMW-trace}
\Tr{2} K_1 = \mu\, I_1, \qquad
\tr_{\! R}\, I ={ \textstyle{ (q-\mu)(q^{-1}+\mu)\over  q-q^{-1}}}.
\ee

Let $\{R,F\}$ be a compatible pair of R-matrices, where
$R$ is skew-invertible of the BMW type and $F$ is strict skew-invertible.
In  \cite{OP} we associated with such a pair an invertible operator $G\in {\rm Aut}(V)$ and 
two invertible linear maps, $\phi$ and $\xi$, acting on the space ${\rm End}(V)\otimes W$ 
where $W$ is an 
arbitrary vector space. We extensively use the operator $G$ and maps $\phi$ and $\xi$ 
in investigations of the BMW type quantum matrix algebras
(see, e.g., sections \ref{secqma}, \ref{secCHSp} below). Here we present formulas for them and for 
their inverses. The operator $G$ and its inverse read
\be
\lb{G}
G_1 \, :=\, \tr_{\!(23)} K_2 F_1^{-1} F_{2}^{-1}, \qquad
G_1^{-1}\, =\, \tr_{\!(23)} F_2 F_1 K_2 .
\ee
The maps $\phi$ and $\xi$ are defined by
\ba
\lb{phi} 
\phi(M)_1 &:=& \Tr{2} \left( F_{1}M_1 F^{-1}_{1} R_{1}\right),
\\[2pt]
\lb{xi}
\xi(M)_1 &:=& \Tr{2} \left( F_{1}
M_1 F^{-1}_{1} K_{1}\right).
\ea
Here $M$ is an arbitrary operator with values in a vector space $W$, $M\in {\rm End}(V)\otimes W$. The inverse maps read
\ba
\lb{phi-inv}
\phi^{-1}(M)_1 &=& \mu^{-2}\TR{2}{R_f} \left(  F_{1}^{-1}
M_1 R^{-1}_{1} F_{1}\right)
\\[2pt]
\lb{xi-inv}
\xi^{-1}(M)_1 &=& \mu^{-2} \TR{2}{R_f}\left( F^{-1}_{1} M_1
K_{1} F_{1}\right) .
\ea
Here the matrix $D_{R_f}$ which is needed for calculations of the $R_f$-traces is
$$
D_{R_f} = D_{F^{-1}} (D_{R^{-1}})^{-1} D_F.
$$

\subsection{Orthogonal and symplectic type R-matrices}\lb{subsec3.4height}

Consider R-matrix realizations $\rho_R(a^{(i)})$ of the antisymmetrizers  (\ref{a^k}).
We impose  additional constraints on a skew invertible BMW-type R-matrix $R$ demanding that
\be
\lb{spec4} 
{\rm rk}\,\rho_R(a^{(i)})\neq 0\, \quad \forall\; i=2,3,\dots ,k\quad\mbox{and} \quad
\rho_R\left( a^{(k)}\sigma^-_{k}(q^{-2k})a^{(k)} \right)\equiv 0\,
\ee
for some $k\geq 2$. Here we assume that the parameters $q$, $\mu$
fulfill conditions (c.f. with the conditions (\ref{mu})$\,$)
\be
\lb{mu1} 
i_q\neq 0\,  \;\;\forall\; i=2,3,\dots ,k;
\qquad \mu\neq -q^{-2i+1}\,  \;\;\forall\; i=1,2,\dots ,k .
\ee
Note that in case $(k+1)_q\neq 0$ the last condition in eq. (\ref{spec4}) means vanishing of the $(k+1)$-st 
antisymmetrizer: $\rho_R(a^{(k+1)})=0$. We do not use this short form to avoid unnecessary restrictions on 
the parameter $q$.

\medskip
An R-matrix satisfying the conditions (\ref{spec4}) is called an
{\em R-matrix of finite height};  the number $k$ is called the
{\em height} of the R-matrix.

\medskip
Let us discuss some consequences of the relations (\ref{spec4}). Applying ${\Tr{i}}$ to $\rho_R(a^{(i)})$
and using the relations (\ref{a^k}), (\ref{Rtr}) and (\ref{BMW-trace}), we calculate
\be\lb{spec1}\Tr{i} \rho_R(a^{(i)})\, =\, \delta_i\, \rho_R(a^{(i-1)})\, \ee
where $\displaystyle{
\delta_i\equiv\delta_i(q,\mu) :=\, - {q^{i-1}(\mu + q^{1-2i})(\mu^2 - q^{4-2i})\over
(\mu + q^{3-2i})(q-q^{-1}) i_q }. }$
In view of eqs.(\ref{spec1}), the last condition in (\ref{spec4}) implies, in particular, that $\delta_{k+1} = 0$,
wherefrom one specifies three admissible values of $\mu$: $\mu\in\{-q^{-1-2k},\pm q^{1-k}\}$.

\smallskip
Notice that the choice $\mu = -q^{1-k}$ contradicts the conditions (\ref{mu1}) in the case when the
number $k$ is even.
In the case when $k$ is odd, the choices $\mu = -q^{1-k}$ and $\mu = q^{1-k}$ are related by a
substitution $R \mapsto  - R$. On the algebra level, this corresponds to an algebra isomorphism (see \cite{OP}, section 2.2)
$\iota'' :{\cal W}_n(q,\mu)\rightarrow {\cal W}_n(-q,-\mu)$, $ \iota''(\sigma_i)= -\sigma_i$,~ $i=1,\dots ,n-1$. The
antisymmetrizers $a^{(i)}$ are invariant under this map.
Therefore we are left with only
two essentially different choices of the parameter $\mu$: either $\mu=-q^{-1-2k}$ or $\mu=q^{1-k}$.
With these
choices, the consistency of the conditions on $\mu$ in eq.(\ref{mu1}) follows from the conditions
on $q$.
\smallskip

We are now ready to define families of the orthogonal and symplectic R-matrices.

\begin{defin}\lb{definition3.11}
Let $R$ be a skew invertible BMW-type R-matrix. Assume additionally that the {\rm{R}}-matrix $R$ has finite height $k$ for some $k\geq 2$.
This implies, in particular, restrictions on $q$: $i_q\ne 0$ for $i=2,\dots ,k$. Then

\noindent
\hspace{10mm}a)~ $R$ is called a $Sp(2k)$-type R-matrix in the case when $\mu=-q^{-1-2k}$;

\medskip
\noindent
\hspace{10mm}b)~ $R$ is called an $O(k)$-type R-matrix in the case when $\mu=q^{1-k}$
and $\mbox{rk}\,\rho_R(a^{(k)})$ $=1$.\end{defin}

{}For the standard R-matrices related to the quantum groups of the series $Sp_q(2k)$ and $SO_q(k)$ \cite{FRT}, the
conditions a) and b), respectively,  and the relations (\ref{spec4}) are fulfilled. This explains
our terminology.

\medskip
The main subject of this paper is an investigation of the general structure  of the quantum matrix algebras associated with the  R-matrices of symplectic type (see next sections). For illustration purposes in subsection \ref{sec4.3} we consider examples of such algebras related to the standard $Sp(2k)$-type R-matrices.
For reader's convenience  we recall formulas for these particular symplectic R-matrices.

\smallskip
The standard $Sp(2k)$-type R-matrix (see \cite{FRT})  reads
\be
\lb{R-Sp}
R^{\mbox{\tiny (st)}} \! :=\!\! \sum_{i,j=1}^{2k} q^{(\delta_{ij}-\delta_{ij'})} E_{ij}\otimes E_{ji} \! +\!
(q-q^{-1}) \!\!\sum_{1\leq j<i}^{2k}\! \bigl\{\!E_{jj}\otimes
E_{ii} \, -\, q^{(\rho_i-\rho_j)}\epsilon_i \epsilon_j \,E_{i`j}\otimes E_{i j'}\!\bigr\}.
\ee
Here $E_{ij}$ are $2k\times 2k$ matrix units; $\delta_{ij}$ is the Kronecker symbol;
\be
\lb{notat-Sp}
i'= 2k+1-i\,;\quad  \epsilon_i=-\epsilon_{i'}=1\, ; \ 
\rho_i = -\rho_{i'}= (k+1-i)\  \forall i: 1\leq i\leq k.
\ee
The corresponding matrices $K^{\mbox{\tiny (st)}}$ and $D_{R^{\mbox{\tiny (st)}}}$ are
\be
\lb{K-Sp}
K^{\mbox{\tiny(st)}} \, =\, \sum_{i,j=1}^{2k} q^{-(\rho_{i}+\rho_{j})} \epsilon_i \epsilon_{j'}E_{ij}\otimes E_{i'j'} , \ 
D_{R^{\mbox{\tiny (st)}}}\, =\, \sum_{i=1}^{2k} q^{-(2k+2\rho_i+1)} E_{ii}\, .
\ee

\begin{rem}\lb{remark3.12.1}
	{\rm For the family of symplectic R-matrices, the case $k=1$ is particular: the antisymmetrizer
		$\rho_R(a^{(2)})$ vanishes and the minimal polynomial of $R$ becomes quadratic.
		The R-matrix $R^{\mbox{\tiny (st)}}$, up to normalization
		and reparameterization $q\mapsto q^{1/2}$,
		is of the Hecke type $GL(2)$ (see $Sp(2)$ examples in the subsection \ref{sec4.3}). This is a manifestation
		of the accidental isomorphism $SL(2)\sim Sp(2)$.
		Accidental isomorphisms for quantum groups, corresponding to the standard deformation, are discussed in \cite{JO}.
	}
\end{rem}

\begin{rem}\lb{remark3.12}
	{\rm Functions
		$$
		\Delta^{(i)}(q,\mu):=\Tr{1,2,\dots ,i}\rho_R(a^{(i)}) = \prod_{j=1}^i \delta_j(q,\mu)
		$$
		are, up to an overall factor, particular elements of a set of rational functions $Q_\lambda(\mu^{-1},q)$
		labelled by partitions $\lambda\vdash i$; we have $\Delta(q,\mu)=\mu^i Q_{[1^i]}(\mu^{-1},q)$. The functions
		$Q_\lambda(\mu^{-1},q)$ were introduced in Theorem 5.5 in \cite{W}.
		They describe the q-dimensions of the
		highest weight modules $V_{\lambda}$ for the orthogonal and symplectic quantum groups (see \cite{W},   
		Section 5
		and \cite{OW}, Lemma 3.1). 
}\end{rem}

\section{Quantum matrix algebra}\lb{secqma}

In this section we recall definitions and main facts about the quantum matrix algebras from \cite{OP}.
A special attention is paid to the family of BMW type quantum matrix algebras. 
The notion of the characteristic subalgebra is introduced and two of its generating sets are described. 
The $\star$-product of the quantum matrices is defined. It substitutes for the usual matrix multiplication in the case of quantum matrices. All these data are necessary for a proper generalization of the Cayley-Hamilton theorem to the case of quantum matrix algebras. The latter is done in the next section.
\medskip

Let $\{R,F\}$ be a compatible pair of R-matrices. In the sequel we assume that $R$ and $F$ are strict skew invertible although some definitions can be given without this condition. 
A {\em quantum matrix algebra} ${\cal M}(R,F)$ is a quotient algebra of the free associative unital algebra $W={\Bbb C}\langle M_a^b\rangle$ by a two-sided ideal generated by entries of the matrix relation
\be 
R_1 M_{\overline 1}M_{\overline 2} = M_{\overline 1}M_{\overline 2}R_1\ . 
\label{qma}
\ee
Here $M = \|M_a^b\|_{a,b=1}^{\mbox{\footnotesize\sc n}}$ is the matrix of generators;
the matrix copies $M_{\overline i}$ are constructed with the help
of the R-matrix $F$ in the following way 
\be
M_{\overline 1}:=M_1, \quad M_{\overline{i}}:=
F^{\phantom{-1}}_{i-1}M_{\overline{i-1}}F_{i-1}^{-1}\ .
\lb{kopii}
\ee
The set of relations 
\be 
R_i M_{\overline i}M_{\overline{i+1}} = M_{\overline i}M_{\overline{i+1}}R_i 
\label{qmai}
\ee
for any given value of the index $i\geq 1$ is equivalent to (\ref{qma}) and can be as well used for the definition of 
the quantum matrix algebra.

\medskip
	Denote by ${\cal C}(R,F)$ a vector subspace of the quantum matrix algebra ${\cal M}(R,F)$
	spanned linearly by the unity and elements
	\be
	\lb{char}
	ch(\alpha^{(n)}) := \Tr{1,\dots ,n}(M_{\overline 1}\dots M_{\overline n}\,
	\rho_R(\alpha^{(n)}))\ ,\quad n =1,2,\dots\ ,
	\ee
	where $\alpha^{(n)}$ is an arbitrary element
	of the braid group ${\cal B}_n$. The space ${\cal C}(R,F)$ is a commutative subalgebra in ${\cal M}(R,F)$ (this is proved in the article \cite{IOP1} which deals with the Hecke type quantum matrix algebras but the proof is valid for an arbitrary compatible 
pair $\{R,F\}$). The algebra ${\cal C}(R,F)$ is called the {\em characteristic subalgebra} of ${\cal M}(R,F)$.

\medskip
Denote by ${\cal P}(R,F)$ a linear subspace of ${\rm End}(V)\otimes {\cal M}(R,F)$
spanned by ${\cal C}(R,F)$-multiples of the identity matrix, $I\, ch\;\; \forall\, ch\in {\cal C}(R,F)$,
and by elements
\be
\lb{pow}
M^1 :=\! M, \ (M^{\alpha^{(n)}})_{1} :=\! \Tr{2,\dots ,n}( M_{\overline 1}
\dots M_{\overline n}\,\rho_R(\alpha^{(n)})),\ n =2,3,\dots,
\ee
where $\alpha^{(n)}$ belongs to
the braid group ${\cal B}_n$. The space ${\cal P}(R,F)$
carries a structure of a right ${\cal C}(R,F)$--module
\be
\lb{r-module}
M^{\alpha^{(n)}} ch(\beta^{(i)}) =M^{(\alpha^{(n)}\beta^{(i)\uparrow n})}\, \ 
\forall\, \alpha^{(n)}\in {\cal B}_n, \; \beta^{(i)}\in {\cal B}_i\, ,\ n,i=1,2,\dots\, ,
\ee
Here, in the right hand side, we denoted by the same symbol $\alpha^{(n)}$ the image of the element $\alpha^{(n)}$
under the natural monomorphism
${\cal B}_{n}\hookrightarrow {\cal B}_{n+i}: \sigma_j\mapsto \sigma_{j}$. The symbol  
$\beta^{(i)\uparrow n}$ in the right hand side denotes the image of the element $\beta^{(i)}$ under the natural monomorphism ${\cal B}_{i}\hookrightarrow {\cal B}_{n+i}: \sigma_j\mapsto \sigma_{j+n-1}$.
Formula (\ref{r-module}) is just a component-wise multiplication of the matrix $M^{\alpha^{(n)}}$ by the element
$ch(\beta^{(i)})$.
\medskip

We call {\it $\star$-product}
the binary operation 
${\cal P}(R,F) \otimes {\cal P}(R,F)\rightarrow \hspace{-4mm}^\star\hspace{3mm} {\cal P}(R,F)$ defined by 
\ba
\nn
( ch(\beta^{(i)}) I) \star  M^{\alpha^{(n)}} :=
M^{\alpha^{(n)}}  ch(\beta^{(i)})&=:&
M^{\alpha^{(n)}}\! \star  ( ch(\beta^{(i)})I)  ,
\\[2pt]
\nn
( ch(\alpha^{(n)})I)\star ( ch(\beta^{(i)})I) &:=& ( ch(\alpha^{(n)}) ch(\beta^{(i)}))I ,
\\[2pt]
\lb{MaMb}
M^{\alpha^{(n)}}\! \star  M^{\beta^{(i)}} 
&:=&  M^{(\alpha^{(n)}\star \beta^{(i)})} ,
\\[2pt]
\nn
\mbox{where we use the notation}\quad\;\;
\alpha^{(n)}\star \beta^{(i)} &:=&
\alpha^{(n)}\beta^{(i)\uparrow n} (\sigma_n\dots \sigma_2 \sigma_1\sigma_2^{-1}\dots \sigma_n^{-1}) .
\ea
The $\star$-product on ${\cal P}(R,F)$ is  associative \cite{OP}. 

\vskip .1cm
In what follows we often use the $\star$-multiplication by the matrix of generators of the quantum matrix algebra  ${\cal M}(R,F)$. Explicitly it reads, see \cite{OP},
\be
\lb{M*}
M \star  N = M\cdot \phi(N)  \quad \forall N\in {\cal P}(R,F),
\ee
where $\cdot$ denotes the usual matrix multiplication and the map $\phi$ is defined in (\ref{phi}). In particular, one can introduce 
the noncommutative analogue of the matrix power:  
\be
\lb{M^k}
M^{\overline{0}} := I\, , \qquad M^{\overline{n}}\, :=\,
\underbrace{M\star  M\star \dots \star  M}_{\mbox{\small $n$ times}}\, =\, M^{(\sigma_1\sigma_2\dots \sigma_{n-1})}.
\ee
Here we use symbol $M^{\overline{n}}$ for the {\em $n$-th power of the matrix $M$}. 

\subsection{BMW type}

If $R$ is an R-matrix of the BMW, $Sp(2k)$ or $O(k)$ type then ${\cal M}(R,F)$ is
called, respectively, a {\em BMW, $Sp(2k)$ or $O(k)$ type quantum matrix algebra}. 

\vskip .1cm
For the BMW type quantum matrix algebra the following relations are satisfied as a consequence of (\ref{qma})
\be
\lb{tau2}
K_{i}\, M_{\overline{i}}M_{\overline{i+1}}\!\! =\!\!
\mu^{-2} K_{i}\, g\, 
\, =\,M_{\overline{i}}M_{\overline{i+1}}\, K_{i}\quad 
\forall\; i\geq 1 , 
\ee
where
\be
\lb{tau}
\displaystyle{ g\! :=\! { { \mu(q-q^{-1})\over (q-\mu)(q^{-1}+\mu) }}\,
\Tr{1,2} \left( M_{\overline{1}}M_{\overline{2}}\, K_1\right) .}
\ee
The element $g$ is called  a   {\em 2-contraction of $M$}. 
	
\vskip .1cm	
For the quantum matrix algebra of the BMW type the 2-contraction $g$ 
is an element of the characteristic subalgebra. The characteristic subalgebra of the BMW type quantum matrix algebra
is generated by either one of the sets  $\{g,p_i\}_{i\geq 0}$, where
\be
\lb{P-i} 
p_0 = \tr_{\!\! R} I ={ \textstyle{ (q-\mu)(q^{-1}+\mu)\over  q-q^{-1}}}  ,\ 
p_1 \!=\! \tr_{\!\! R}\, M,
\ p_i =  ch(\sigma_{i-1}\dots\sigma_2\sigma_1) \ i=2,3,\ldots\ ,
\ee
or $\{g,a_i\}_{i\geq 0}$, where \ba
\lb{A_i}
a_0 \!&=&\! 1,\qquad\quad\;\, a_i = ch(a^{(i)}) \qquad\qquad\; i=1,2,\dots  .
\ea
Elements $p_i$ and $a_i$ are called {\em power sums} and {\em elementary symmetric functions}, respectively.

For the Hecke type quantum matrix algebra, the corresponding algebra  ${\cal P}(R,F)$, as the ${\cal C}(R,F)$-module, is spanned by the matrix powers $M^{\overline n}$, $n\geq 0$, of the generating matrix $M$. For the BMW type quantum matrix algebra this is not the case.  Namely, as the ${\cal C}(R,F)$--module,  the BMW type algebra ${\cal P}(R,F)$ is spanned by matrices
(see \cite{OP}, proposition 4.11)
$$
M^{\overline{n}}\,  \quad \mbox{and}\quad
\Mt(M^{\overline{n+2}})\, , \quad n=0,1,\dots\, 
$$
Here we introduced a ${\cal C}(R,F)$--module map~
$\Mt : {\cal P}(R,F)$ $\rightarrow$ $ {\cal P}(R,F)$
\be
\lb{Mt}
\Mt (N) := M\cdot \xi(N), \qquad  \forall\, N\in {\cal P}(R,F),
\ee
where the map $\xi$ is given in (\ref{xi}).\smallskip
 
The BMW type algebra ${\cal P}(R,F)$ is commutative  \cite{OP}. \medskip

\vspace{0mm}

To define inverse powers of the quantum matrix $M$ one considers the extension of the BMW type algebra 
${\cal M}(R,F)$ by the inverse  $g^{-1}$ of the 2-contraction 
\be
\lb{j-inv} 
g^{-1}\, g\, =\, g\, g^{-1}\, =\, 1\, , \qquad g^{-1}\, M\, =\, (G^{-1}MG)\, g^{-1}\, .
\ee
where the numeric matrices $G^{\pm 1}\in {\rm Aut}(V)$ are defined in eqs. (\ref{G}). The latter relation in (\ref{j-inv}) is justified by the permutation rules for the 2-contraction. For an arbitrary matrix
$N\in {\cal P}(R,F)$ it reads
\be
\lb{g-perm}
N\, g\, =\, g\, ( G^{-1}N G) .
\ee
\noindent {\bf Proof.}~ In a particular case $N=M$ --- the matrix of generators of ${\cal M}(R,F)$ this formula is proved in \cite{OP}, lemma 4.13. Consequently, by lemma 3.11, eq. (3.45), \cite{OP}, we have $M_{\overline j}\,g\, =\, g\, ( G_j^{-1}M_{\overline j}\, G_j)$, $j=1,2,\ldots$ Thus for $u=M_{\overline 1}\dots M_{\overline n}\,\rho_R(\alpha^{(n)}))$, $\alpha^{(n)}\in {\cal B}_n$, we have
$ug=g G_1G_2\dots G_n u G_n^{-1}\dots G_2^{-1}G_1^{-1}$. By the cyclic property of the trace and lemma 3.11, 
eq. (3.44),  $G_2\dots G_n$ cancels with $G_n^{-1}\dots G_2^{-1}$ which proves eq.(\ref{g-perm}) for $N=\Tr{2,\dots ,n}(u)$.\hfill$\blacksquare$

\vskip .1cm
The extended algebra, which we shall further denote by ${\cal M^{^\bullet\!}}(R,F)$, contains
the inverse matrix to the matrix $M$
\be
\lb{M-inv} 
M^{-1}\, =\, \mu\, \xi(M)\, g^{-1},\qquad
M\cdot M^{-1}\,=\, I\, =\, M^{-1}\cdot M.
\ee 
The matrix $M^{-1}$ is the inversion of $M$ with respect to the usual matrix product. Inversion with respect 
to the $\star$-product looks differently 
\be
\lb{M-star-inv}
M^{\overline{-1}} =
\phi^{-1}(M^{-1}),\qquad M^{\overline{-1}}\star M \, =\, I\,  =\, M\star M^{\overline{-1}}.
\ee
In general, $M^{\overline{-1}}\neq M^{-1}$.\medskip

One can define the unique  
extension ${\cal P^{^\bullet\!}}(R,F)$ of the algebra ${\cal P}(R,F)$
by a repeated $\star$-multiplication with
$M^{\overline{-1}}$
\be
\lb{Minv*N}
M^{\overline{-1}}\star  N \,:=\, \phi^{-1}(M^{-1}\cdot N)\, =:\, N \star  M^{\overline{-1}}
\qquad
\forall\; N\in{\cal P^{^\bullet\!}}(R,F)\, .
\ee
The algebra ${\cal P^{^\bullet\!}}(R,F)$ is associative and commutative with respect to the $\star$-product. It is 
also the right ${\cal C^{^\bullet\!}}(R,F)$-module algebra with respect to the extension 
${\cal C^{^\bullet\!}}(R,F)\supset {\cal C}(R,F)$ of the characteristic subalgebra by the element $g^{-1}$.

Particular examples of the $\star$-multiplication by $M^{\overline{-1}}$ are the inverse $\star$-powers
of $M$
$$
M^{\overline{-n}}\, :=\,\underbrace{M^{\overline{-1}}\star \dots \star  M^{\overline{-1}}\star }_{n\ \text{times}}I. $$
The $\star$-powers obey the usual rules of the $\star$-product of matrix powers:
$M^{\overline{i}}\star M^{\overline{n}}\, =\, M^{\overline{i+n}}\;\; \forall\  i,n\in {\Bbb Z}$.

\section{Cayley-Hamilton theorem}\lb{secCHSp}

The Cayley-Hamilton theorem for the orthogonal and symplectic quantum groups was stated in the unpublished text \cite{OP2}. Here 
we establish and discuss in details a strengthened version of the Cayley-Hamilton theorem in the symplectic case.

\vskip .1cm
Throughout this section we assume that $\{R,F\}$ is a compatible pair of R-matrices, in which the operator
$F$ is strict skew invertible and the operator $R$ is skew invertible of the BMW-type and,
hence, strict skew invertible.

\vskip .1cm
In the subsection \ref{subsec5.1} we investigate matrix relations in the algebra ${\cal P}(R,F)$ involving `wedge' powers of the quantum matrix $M$: $M^{a^{(i)}}$, $0\leq i\leq n$. We confine the eigenvalues $q$ and $\mu$ of the matrix $R$ by conditions
$$
i_q\neq 0,\;\mu\neq -q^{3-2i}\;\;\forall
\; i=2,3,\dots ,n,
$$
in which case all the antisymmetrizers $a^{(i)}\in {\cal W}_{n}(q,\mu)$, $i=2,3,\dots ,n$,
 and, hence, the elements
$a_i\in {\cal C}(R,F)$ and the matrices $M^{a^{(i)}}\in {\cal P}(R,F)$  are well defined.

\vskip .1cm
Conditions on $R$ specific for the R-matrices of the type $Sp(2k)$,  are imposed in Subsection \ref{subsec5.2}.

\subsection{Basic identities}\lb{subsec5.1}

Consider a set of `wedge' powers of the quantum matrix $M$: $M^{a^{(i)}}\in {\cal P}(R,F)$. Following \cite{OP},
we introduce series of matrices in ${\cal P}(R,F)$, which we further refer to as `descendants' of the matrices $M^{a^{(i)}}$.
\be
\lb{La}
\begin{array}{l}
A^{(m,i)}\ :=\ i_q\, M^{\overline{m}}  \star  M^{a^{(i)}}  \\[1em] 
B^{(m+1,i)}\ :=\
i_q\, M^{\overline{m}} \star  \Mt (M^{a^{(i)}})\end{array}
\quad   \forall\, i,m:\; 1\leq i\leq n, \;m\geq 0.
\ee
It is  suitable to set, by definition,
\ba
\lb{LT0}
&A^{(m,0)}\ :=0\ \ \ {\mathrm{and}}\ \ \ \ B^{(m,0)}\ :=\ 0\,\qquad \forall\; m\geq 0\, .&\ea
and to complement the series by the elements\footnote{
Note that $A^{(-1,i)}$ and $B^{(0,i)}$ belong to the extension of the algebra ${\cal P}(R,F)$ by the
$\star$-inverse matrix $M^{\overline{-1}}$.
}
\be
\lb{AB-boundary}
\begin{array}{l}
A^{(-1,i)}\, :=\, i_q\, \phi^{-1}\left(
\Tr{2,3,\dots i} M_{\overline{2}}M_{\overline{3}}\dots M_{\overline{i}}\, \rho_R(a^{(i)})
\right)\ , \\[1em]
B^{(0,i)}\ :=\ i_q\, \phi^{-1}\bigl(\xi\bigl(M^{a^{(i)}}\bigr)\bigr)\, . \end{array}
\ee
The following recursive relations among the descendants are derived in \cite{OP}:
\begin{lem}\lb{lemma5.1}
{}For ~$0\leq i\leq n-1$~ and ~$m\geq 0$, the matrices $A^{(m-1,i+1)}$ and $B^{(m+1,i+1)}$
satisfy equalities
\ba
\lb{rek1}
A^{(m-1,i+1)} &=& q^i M^{\overline{m}}\, a_i\, -\,
A^{(m,i)}\, -\, { \mu q^{2i-1}(q-q^{-1})\over 1+\mu q^{2i-1}}\ B^{(m,i)}\, ,
\\[1em]
\lb{rek2}
B^{(m+1,i+1)} &=&\Bigl( \mu^{-1}q^{-i} M^{\overline{m}}\, a_i\, +\,
{q-q^{-1}\over 1+\mu q^{2i-1}}\ A^{(m,i)}\, -\,B^{(m,i)}\Bigr) g\,  .
\ea
\end{lem}

\medskip
By a repeated use of these recurrent relations one can derive for a certain subset of the descendants
their expansions in terms of non-negative matrix powers $M^{\overline{j}}$, $j\geq 0$, only\footnote{By
Proposition 4.11 \cite{OP}, one expects also presence of the terms $\Mt(M^{\overline{j}})$
in the expansions of generic descendants.}. For the Hecke type QM-algebras analogues of these expansions are known as the Cayley-Hamilton-Newton identities
\cite{IOP,IOP1,IOPS}.

\begin{prop}\lb{corollary5.2}
{}For $1\leq i\leq n$ and $m\geq i-2$, one has
\be
\lb{cor1a}
A^{(m,i)}\;\, =\;\, (-1)^{i-1}\sum_{j=0}^{i-1} (-q)^j\Bigl\{
M^{\overline{m+i-j}} + {1-q^{-2}\over 1+\mu q^{2i-3}}\sum_{r=1}^{i-j-1}M^{\overline{m+i-j-2r}} (q^2 g)^r
\Bigr\} a_j\, .
\ee
\hspace{21mm}For $1\leq i\leq n$ and $m\geq i$, one has
\ba
\nonumber
B^{(m,i)}&=& (-1)^{i-1}\sum_{j=0}^{i-1} (-q)^j\Bigl\{
\mu^{-1} q^{-2j} M^{\overline{m-i+j}} g^{i-j}\hspace{56mm}
\\[0em]
&&\hspace{45mm}
-{q^{-1}(1-q^{-2})\over 1+\mu q^{2i-3}}\sum_{r=1}^{i-j-1}M^{\overline{m+i-j-2r}} (q^2 g)^r
\Bigr\} a_j\, .
\lb{cor1b}
\ea
\end{prop}

\noindent {\bf Proof.}~ We employ induction on $i$.
In the case $i=1$, the relations (\ref{cor1a}) and (\ref{cor1b}) reproduce the definitions (\ref{La}):
$$
A^{(m,1)}\, =\, M^{\overline{m+1}}\ , \qquad B^{(m,1)}\, =\, \mu^{-1} M^{\overline{m-1}} g\ .
$$

\smallskip
It is then straightforward to verify the induction step $i\rightarrow i+1$ with the help of the relations
(\ref{rek1}) and (\ref{rek2}). \hfill$\blacksquare$

\begin{rem}
{\rm
When ~$m\geq i-2$~ (respectively, ~$m\geq i$),
all the $\star\, $-powers of $M$ in the right hand side of the relation (\ref{cor1a})
(respectively, the relation (\ref{cor1b})$\,$) are non-negative. This is why we specify these restrictions on $m$. For
an invertible matrix $M$, the restrictions on $m$ can be removed.
}
\end{rem}

\begin{rem}
{\rm
The Hecke type version of these relations can be reproduced by setting $g=0$ in formulas of Proposition \ref{corollary5.2}.
Relation for $B^{(m,j)}$ becomes trivial. Relation (\ref{cor1a}) for $A^{(m,j)}$ simplifies drastically, the terms
with the element $g$ disappear and the condition $m\geq i-2$ weakens to $m\geq -1$.
{}For $m=0$, the relation (\ref{cor1a}) reproduces the Cayley--Hamilton--Newton identities found in \cite{IOP,IOP1}.
The R-trace maps
of these identities are the Newton relations.
In the $GL(k)$-case, that is, if the operator $R$ fulfills the condition $\rho_R(a^{(k+1)})=0$, the left hand side of the relation (\ref{cor1a}) vanishes
in the case $i=k+1$. Then, with the choice $m=-1$ the relation (\ref{cor1a})
reproduces the Cayley--Hamilton identity.}
\end{rem}

\subsection{Cayley-Hamilton theorem: type $Sp(2k)$}\lb{subsec5.2}

Specifying to the case of the $Sp(2k)$-type quantum matrix algebra, we notice that the 
condition $\mu= - q^{-1-2k}$ leads to the following linear dependency
between $A^{(m-1,k+1)}$, see (\ref{rek1}), and  $B^{(m+1,k+1)}$, see (\ref{rek2}):
\be\lb{previous-rem}
\left. \left(B^{(m+1,k+1)} + q A^{(m-1,k+1)} g \right)\right\vert_{\mu=-q^{-1-2k}} = 0\,\ \quad \forall\; m\geq 0\, .
\ee
The height $k$ condition (\ref{spec4}) on the $Sp(2k)$-type R-matrix $R$ cuts the series of 'descendants' $A^{(m,i)}$ and $B^{(m,i)}$
at the level $i=k+1$: $A^{(m-1,k+1)} = B^{(m+1,k+1)}=0\, \quad \forall\; m\geq 0$. 
The Cayley-Hamilton theorem follows exactly from these cutting conditions. The relations  (\ref{previous-rem}) show that all the conditions for $B^{(m+1,k+1)}$ follow from the conditions for $A^{(m-1,k+1)}$. In turn, by eqs. (\ref{La}) and (\ref{AB-boundary}) we have 
\be
\lb{ahah-new}
A^{(m-1,k+1)}=M^{\overline{m}}  \star A^{(-1,k+1)}.
\ee 
Thus, all the cutting conditions arise from the single one 
\be
\lb{ahah} A^{(-1,k+1)}\, =\, 0\, .
\ee
Unfortunately, the latter condition cannot be expressed in terms of nonnegative powers of the matrix $M$ only. By Proposition \ref{corollary5.2}, for the condition
\be
\lb{ahah-1}
A^{(k-1,k+1)}\, =\,0
\ee
such an expression does exist.

The relations (\ref{ahah}) and (\ref{ahah-1}) are equivalent if the 2-contraction $g$ and, hence, the matrix $M$ are invertible.
We shall first investigate the condition (\ref{ahah-1}).
Substituting $\mu=-q^{-1-2k}$ and (\ref{cor1a}) into (\ref{ahah-1})
and rearranging terms of the sum  we obtain the Cayley-Hamilton theorem for the quantum matrices of the type $Sp(2k)$:

\begin{theor}
\lb{theorem5.4}
Let ${\cal M}(R,F)$ be the $Sp(2k)$-type quantum matrix algebra.
Then the quantum matrix $M$  of the algebra generators satisfies the Cayley-Hamilton identity
\be
\lb{CHSp-1}
\sum_{i=0}^{2k} (-q)^i M^{\overline{2k-i}} \epsilon_i\, =\,  0\, ,
\ee
where
\be
\lb{CHSp-2}
\epsilon_i\, :=\, \sum_{j=0}^{[i/2]} a_{i-2j}\, g^j\, , \qquad
\epsilon_{k+i}\, :=\, \epsilon_{k-i}\, g^{i}\,\qquad \forall\; i=1,2,\dots ,k\, .
\ee
\end{theor}

Let us now consider the matrix identity (\ref{ahah}). In case of non-invertible $g$, this identity is more informative than the Cayley-Hamilton identity (\ref{CHSp-1}).
Matrix components of its left hand side are $k$-th order homogeneous polynomials in the components of the quantum matrix $M$, containing, apart of $M$, $\star$-powers of yet another quantum matrix obtained from $M$ by a linear map $\pi := \mu\, \phi^{-1}\!\circ \xi$ (see eqs. (\ref{xi}), (\ref{phi-inv})).

\begin{lem}
\lb{lem4.6}
For the compatible pair $\{R,F\}$ of strict skew invertible R-matrices, where $R$ is of the BMW-type, the map $\pi := \mu\, \phi^{-1}\!\circ \xi$ does not depend on $F$.
The explicit formulas for $\pi$  and  $\pi^{-1}$  read:
\ba
\lb{pi}
\pi(M)_1 &=&  \Tr{2} R_{12} M_1 K_{12}\, =\, \Tr{2} K_{12} M_1 R_{12}\, .
\\[2pt]
\lb{pi-inv}
\pi^{-1}(M)_1&=& \mu^{-2}\, \Tr{2} R_{12}^{-1} M_1 K_{12}\, =\, \mu^{-2}\, \Tr{2} K_{12} M_1 R^{-1}_{12} .
\ea
\end{lem}

\noindent {\bf Proof.}~ Instead of proving the first equality in (\ref{pi}) directly it is easier to verify the relation $\phi (\Tr{2} R_{1} M_1 K_{1})=\mu\, \xi(M)_1$:
\ba
\nn
\phi\bigl( \Tr{2} R_{1} M_1 K_{1} \bigr)
& =&
 \Tr{2} F_{12} \Bigl\{ \Tr{2'} R_{12'} M_1 K_{12'} \Bigr\} F^{-1}_{12} R_{12}
\\[2pt] 
\nn
&=&\, \Tr{23} \underline{F_1\bigl\{ F_2 R_1} M_1 K_1 \underline{F_{2}^{-1}\bigr\} F_1^{-1}} R_1
\\[2pt]
\nn
&=&
 \Tr{23} \underline{R_2}  F_1 \underline{F_2 M_1 F_{2}^{-1}} F_1^{-1}  K_2 R_1 \\[2pt]
 &=&\, \Tr{23}   F_1  M_1  F_1^{-1}  \underline{K_2 R_1 R_2}
\\[2pt]
\nn
&= &
\Tr{2\underline{3}}   F_1  M_1  F_1^{-1}  \underline{K_2} K_1\\[2pt]
\nn
&=&\, \mu \Tr{2}   F_1  M_1  F_1^{-1}  K_1 \, =\, \mu\, \xi(M).
\ea
Here in calculations we underline terms which undergo a transformation in the next step. For the transformations
we used the compatibility relations for the pair $\{R,F\}$ (\ref{sovm}), BMW algebra relations for the matrices $R$ and $K$ (\ref{bmwRa}), first formula in (\ref{BMW-trace}), and the following properties of the R-trace (see \cite{OP}, lemma 3.2 and corollary 3.4)
\[
\!\!\!\!\Tr{2} F_1^{\pm 1}\, X_1\, F_1^{\mp 1}= I_1\, \tr_{\! R} X 
\quad\; \forall\; X\in {\rm End}(V)\otimes W, \]
 where $W$ is a $\Bbb C$-linear space, and
\[
\Bigl[ R_{12}, \, {D_R}_1 {D_R}_2 \Bigr]  = 0
\]
which is equivalent to 
\[ \Tr{12} Y_1\, R_1 = \Tr{12} R_1\, Y_1
\qquad \forall\; Y\in {\rm End}(V^{\otimes 2})\otimes W.\]
The second equality in relation (\ref{pi}) holds for an arbitrary BMW type R-matrix  $R$. 

\vskip .1cm
To prove formula (\ref{pi-inv}) one notices that the map $\pi$ is proportional to the map $\xi$ for the pair $\{R,R\}$. Thus the first equality in (\ref{pi-inv}) follows from the formula
for  $\xi^{-1}$ for the pair $\{R,R\}$, see (\ref{xi-inv}). The second equality  in (\ref{pi-inv}) is obtained from the first one by the same remark as for the map $\pi$. 
\hfill$\blacksquare$
\medskip

Until the end of this subsection we let $M$ to be the matrix of generators of the BMW type quantum matrix algebra ${\cal M}(R,F)$.

\vskip .1cm
In general, the matrix $\pi(M)$ does not belong to the  algebra ${\cal P}(R,F)$.
On the other hand,  $\pi(M)$ is related to the $\star$-inverse of the matrix $M$ (see (\ref{M-star-inv}))
\be
\lb{pi-Minv}
\pi(M)\, =\, M^{\overline{-1}} g\, =\, M^{\overline{-1}}\star I g
\ee
and thus belongs to the extended algebra ${\cal P}^{^\bullet}\!(R,F)$.
 The formula for the $\star$-product for the matrix $\pi(M)$ is clearly induced from that for $M^{\overline{-1}}$ (see (\ref{Minv*N})) and the permutation rules for $g$ (see (\ref{g-perm})):
\be
\lb{piM-star}
\pi(M)\star N\, := N\star \pi(M)\, :=\, 
\mu \phi^{-1}( \xi(M)\cdot G^{-1} N G), \qquad \forall\; N\in {\cal P}(R,F).
\ee
Complementing the  algebra ${\cal P}(R,F)$ with the $\star$-multiples of $\pi(M)$
$$
\pi(M)^{\overline{n}}\; :=\;  \underbrace{\pi(M)\star \dots \star  \pi(M) }_{n\ times}
$$
one obtains an intermediate extension ${\cal P}^\circ(R,F)\supset {\cal P} (R,F)$,
${\cal P}^\circ (R,F)\subset  {\cal P}^\bullet (R,F)$. It is this algebra where the matrix $A^{(-1,k+1)}$ belongs to.

Now we are ready to  write down the identity (\ref{ahah})  in terms of $\star$-powers of the matrices $M$ and $\pi(M)$.
\begin{prop}
\lb{preCH}
Let ${\cal M}(R,F)$ be the $Sp(2k)$-type quantum matrix algebra. The matrix $M$ of generators of this algebra and its image $\pi(M)$ under the map (\ref{pi})
satisfy the following $k$-th order matrix polynomial identity
\be
\lb{preCHSp}
\sum_{i=0}^{k} (-q)^i M^{\overline{k-i}} \epsilon_i\, + \, q^{2k}\sum_{i=0}^{k-1} (-q)^{-i} \pi(M)^{\overline{k-i}} \epsilon_i\,=\,  0\, .
\ee
Here the coefficients $\epsilon_i$, $i=1, \dots , k,$ are given by eq. (\ref{CHSp-2}).
\end{prop}

\subsection{Simple examples and classical limit.}
\lb{sec4.3}

In this section we present the Cayley-Hamilton and `pre-Cayley-Hamilton' identities (\ref{CHSp-1})  and (\ref{preCHSp}) for the standard RTT- and RE-algebras corresponding to the
$Sp(2k)$-type R-matrix (\ref{R-Sp}) in cases $k=1,2$.
\medskip

{\em Standard $Sp(2)$-type RTT-algebra} is the quantum matrix algebra ${\cal M}(R^{\mbox{\tiny (st)}},P)$, where the R-matrix $R^{\mbox{\tiny (st)}}$ (\ref{R-Sp}) and permutation $P$ act on a tensor square of the 2-dimensional vector space. We use the symbol $T$ for the $2\times 2$ matrix of generators of this algebra.
Permutation relations for its components $T^i_j$, $i,j\in\{1,2\}$ are identical to the permutaton relations of the standard $GL_{q^2}(2)$-type
RTT-algebra:
\be
\lb{RTT-Sp2}
q^2 T^{i}_{2} T^{i}_{1}\, =\, T^{i}_{1} T^{i}_{2}, \quad  q^2 T^{2}_{i} T^{1}_{i}\, =\, T^{1}_{i} T^{2}_{i},\quad [T^{2}_{1}, T^{1}_{2}] \, =\, 0,\quad
[T^{2}_{2}, T^{1}_{1}] \, =\, (q^{-2}-q^2) T^{1}_{2} T^{2}_{1}.
\ee 
The R-matrix image $\rho_{{R^{\mbox{\tiny(st)}}}}(a^{(2)})$ of the second order antisymmetrizer vanishes in this 
particular case and, therefore, there is no any additional $g$-covariance conditions.

The two generators of the characteristic subalgebra $g$ and $a_1$ read
\ba
\nn
g& =& {q^{-6}\over q^2+q^{-2}} \left( q^{-2}\, T^1_1 T^2_2 +q^{2}\, T^2_2 T^1_1 - T^1_2 T^2_1 -T^2_1 T^1_2\right)
\\
\lb{g-RTT-Sp2}
 &=&
q^{-6}\left( T^1_1 T^2_2 -q^2\, T^1_2 T^2_1\right),\\[2pt]
\nn
a_1& =& {\rm Tr}_{_{\! R}} T \, =\, q^{-5}\, T^1_1 + q^{-1}\, T^2_2,
\ea
where the second simplified expression for $g$ is obtained with the help of the permutation relations (\ref{RTT-Sp2}).
The 2-contraction $g$ is central in this case (the matrix $G$ for the R-matrix pair $\{R^{\mbox{\tiny (st)}},P\}$  equals the unity), while the element $a_1$ is not.

\vskip .1cm
To write down the characteristic identities for this algebra
we need explicit expressions for the maps $\phi$ (\ref{phi}) and $\pi$ (\ref{pi})
\be
\lb{phi-pi-Sp2}
\phi(T)=
\left(\!\!
\begin{array}{cc}
\scriptstyle
q^{-4}\, T^1_1+{ (1-q^{-4})}\,T^2_2 &
\scriptstyle
q^{-6}\, T^1_2\\[2pt]
\scriptstyle
q^{-2}\, T^2_1 &
\scriptstyle
T^2_2
\end{array}\!\!
\right)\! ,
\quad
\pi(T)=
\left(\!\!
\begin{array}{cc}
\scriptstyle
(q^{-6}+q^{-2})\, T^1_1+q^{-2}\,T^2_2 &
\scriptstyle
-q^{-2}\, T^1_2\\[2pt]
\scriptstyle
-q^{-2}\, T^2_1 &
\scriptstyle
q^{-6}\,T^2_2
\end{array}\!\!
\right)\! .
\ee
The Cayley-Hamilton identity (\ref{CHSp-1}) and its parent identity (\ref{preCHSp}), respectively, read
\ba
\lb{CH-RTT-Sp2}
T^{\overline{2}}\, -\, q\, T\, a_1 \, +\, q^2\, I\, g\,=\, 0,
\\[2pt]
\lb{preCH-Sp2}
T\, -\, q\, I\, a_1 \, +\, q^2\, \pi(T)\, =\, 0,
\ea
where $T^{\overline 2}= T\cdot \phi(T)$.

\vskip .1cm
We note that the identity (\ref{CH-RTT-Sp2}) coincides with the Cayley-Hamilton identity for the standard 
$GL_{q^2}(2)$-type RTT-algebra (see \cite{EOW,IOP,IOP2}), where the 2-contraction $g$ plays the role of the 
quantum determinant of the matrix $T$.
In this particular case the Cayley-Hamilton identity (\ref{CH-RTT-Sp2}) encodes the half of the permutation 
relations (\ref{RTT-Sp2}); in general, a half-quantum  matrix of $GL$ type satisfies the Cayley-Hamilton identity
\cite{CFR,IO}.

\vskip .1cm
Another specific feature of the $Sp(2)$ case is that the `parent' Cayley-Hamilton identity (\ref{preCH-Sp2}) being 
linear in generators is satisfied without any reference to the quadratic permutation relations. 

\medskip
{\em The standard $Sp(2)$-type Reflection Equation (RE) algebra} is the quantum matrix algebra 
${\cal M}(R^{\mbox{\tiny (st)}},R^{\mbox{\tiny (st)}})$, where the R-matrix $R^{\mbox{\tiny (st)}}$ (\ref{R-Sp}) 
acts on the
tensor square of the 2-dimensional vector space. We use the symbol $L$ for the $2\times 2$ 
matrix of generators of 
this algebra. The permutation relations for its components
$L^i_j$, $i,j\in\{1,2\}$, are identical to the permutation relations for the standard $GL_{q^2}(2)$-type RE-algebra:
\ba
\nn
&
L^{i}_{j} L^{1}_{1}\, =\, q^{4(j-i)}\,L^{1}_{1} L^{i}_{j}, \qquad  [ L^{2}_{2}, L^{1}_{2}]\, =\, (1-q^{-4})\,L^1_1 L^1_2,&\\ 
\lb{REA-Sp2}
&
[ L^{2}_{2}, L^{2}_{1}]\, =\, -q^{-4}(1-q^{-4})\,L^1_1 L^2_1,
&
\\[2pt]
\nn
&
[L^{2}_{1}, L^{1}_{2}] \, =\, (1-q^{-4})\, L^1_1(L^{1}_{1} - L^{2}_{2}).
&
\ea 
The two generators of the characteristic subalgebra $g$ and $a_1$ are
\ba
\nn
g& =& {q^{-4}\over q^2+q^{-2}} \left(  L^1_1 L^2_2 + L^2_2 L^1_1 -(1-q^{-4})\,(L^1_1)^2 - L^1_2 L^2_1 -q^4\, L^2_1 L^1_2\right)
\\
\lb{g-REA-Sp2}
 &=&
q^{-2}\left( L^1_1 L^2_2 -(1-q^{-4})\, (L^1_1)^2-L^1_2 L^2_1\right),\\[2pt]
\nn
a_1& =& {\rm Tr}_{_{\! R}} L \, =\, q^{-5}\, L^1_1 + q^{-1}\, L^2_2,
\ea
where the second expression for $g$ is obtained with the use of permutation relations (\ref{REA-Sp2}).
As for any RE-algebra, the generators $g$ and $a_1$ are central.

\vskip .1cm
Another distinguishing property of the RE-algebras --- the identity of the map $\phi$ --- makes their characteristic identity  (\ref{CH-RTT-Sp2}) particularly simple and similar to the classical case. In our situation it reads
\be
\lb{CH-REA-Sp2}
L^2\, -\, q\, L\, a_1 \, +\, q^2\, I\, g\,=\, 0,
\ee
where $L^2$ means the usual matrix square of $L$ and the coefficients $g$ and $a_1$ are given by (\ref{g-REA-Sp2}). Again, this matrix equality encodes half of the permutation relations (\ref{REA-Sp2}).

\vskip .1cm
As stated in lemma \ref{lem4.6} the map $\pi$ depends on the first R-matrix from the compatible pair $\{R,F\}$ only. Hence, for the RTT- and RE-algebra generating matrices $T$ and $L$ the map $\pi$ is literally the same (see (\ref{phi-pi-Sp2})), and the parent Cayley-Hamilton identities for the Sp(2)-type RTT- and RE-algebras coincide (see \ref{preCH-Sp2}).
\medskip

Next, we consider a less trivial example in order to demonstrate the results of this section in 
a greater generality.
It is the
{\em standard $Sp(4)$-type RTT-algebra} --- the quantum matrix algebra ${\cal M}(R^{\mbox{\tiny (st)}},P)$, where 
the R-matrix $R^{\mbox{\tiny (st)}}$ (\ref{R-Sp}) and the permutation $P$ now act on the tensor square 
of the 4-dimensional vector space. We keep notation $M$ for the $4\times 4$ matrix of generators of this algebra. 
Quadratic relations in this algebra consist of 120 permutation relations for 
16 matrix components, and of 10  additional conditions.  
The latter ones together with expression for 
the 2-contraction $g$ can be extracted from the matrix equalities
 (\ref{tau2}), where $i=1$ and
 $\mu=-q^{-5}$ in our case. All the quadratic relations and the expressions for $g$ are collected in the Appendix. There and in the formulas below it is suitable to 
 break  the $4\times 4$ matrix $M$ into four $2\times 2$ blocks $A$, $B$, $C$ and $D$:

\be
\lb{M-part}
M\, =\,
\left(\!\!
\begin{array}{cc}
	A &  B
	\\
	C & D
\end{array}\!\!
\right).
\ee
The coefficients $\epsilon_1$ and $\epsilon_2$ of the Cayley-Hamilton identity, together with the 2-contraction $g$ generate the characteristic subalgebra. The 2-contraction $g$ is central, while $\epsilon_i$, $i=1,2$, are not. Expression for $g$ is given in the Appendix (see eq.(\ref{g-RTT-Sp4})); formulas for $\epsilon_i$ read
\ba
\nn
\epsilon_1 &=& a_1\, =\, q^{-9} A^1_1+q^{-7} A^2_2 +q^{-3} D^1_1 + q^{-1} D^2_2,
\\[3pt]
\nn
\epsilon_2 &=& a_2 + g \, =\, q^{-16}(A^1_1 A^2_2 - q A^1_2 A^2_1) \\
\nn
&&\phantom{ a_2 + g \, =\,} + q^{-4} (D^1_1 D^2_2 -q D^1_2 D^2_1)+ 
q^{-12} (D^1_1 +q^2 D^2_2)(A^1_1+q^2 A^2_2)
\\[2pt]
\nn
&&\phantom{ a_2 + g \, =\,} -  q^{-12} \left(q^{-1}C^1_1 B^1_1 -(q-q^{-1})C^1_1 B^2_2+ C^1_2 B^2_1 +C^2_1 B^1_2 +q^3 C^2_2 B^2_2\right) .
\ea
To write down the characteristic identities we also need expressions for the maps $\xi^{\pm 1}$, $\phi^{\pm 1}$. They are
\ba
\nn
\xi(M)& =& 
\left(\!
\begin{array}{cc}
-q^{-5} \sigma_q(D) &  q^{-8}\sigma_q(B)
\\[2pt]
q^{-2} \sigma_q(C) &-q^{-5} \sigma_q(A)
\end{array}\!
\right),
\\[3pt]
\nn
\phi(M)& =& 
\left(\!
\begin{array}{cc}
q^{-6} \alpha^+_q(A)+(1-q^{-2})\beta_q(D) & q^{-7} \alpha^-_q(B)
\\[2pt]
q^{-1} \alpha^-_q(C) & \alpha^+_q(D)
\end{array}\!
\right),
\\[3pt]
\lb{xi=phi-inv-Sp4}
\xi^{-1}(M)& =& \xi(M)|_{q\leftrightarrow q^{-1}},
\qquad
\phi^{-1}(M)\, =\, \phi(M)|_{q\leftrightarrow q^{-1}}.
\ea
Here $\sigma_q$, $\alpha^\pm_q$, $\beta_q$ are linear maps of the $2\times 2$ matrices
\ba
\nn
\sigma_q(X)=
\left(\!\!
\begin{array}{cc}
\scriptstyle  X^2_2 &\scriptstyle q^{-1} X^1_2
\\
\scriptstyle q X^2_1 &\scriptstyle X^1_1
\end{array}\!\!
\right),
\;\;
\alpha^{\pm}_q(X)=
\left(\!\!
\begin{array}{cc}
\scriptstyle q^{-2} X^1_1\,\pm\, (1-q^{-2}) X^2_2 &\scriptstyle q^{-3} X^1_2
\\
\scriptstyle q^{-1} X^2_1 &\scriptstyle X^2_2
\end{array}\!\!
\right)\!,
\ea
\be
\beta_q(X) = {\scriptstyle (q^{-2} X^1_1+X^2_2)} I + {\scriptstyle q^{-4}}\sigma_q(X),
\ee
The following properties of these maps make the check of the relations (\ref{xi=phi-inv-Sp4}) staightforward:
\ba
\nn
&&(\sigma_q)^{-1}\, =\,\sigma_{1/q},\qquad\quad
(\alpha^\pm_q)^{-1}\, =\, \alpha^\pm_{1/q},
\qquad\quad
\beta_q\circ \alpha^+_{1/q} \,=\,  q^{-4}\, \alpha^+_q \circ \beta_{1/q}.
\ea
The composite map $\pi=-q^{-5} (\phi^{-1}\circ \xi)(M)$  reads explicitly
\be
\nn
\pi(M) = 
\left(\!
\begin{array}{cc}
	q^{-4}( \alpha^+_{1/q}\circ \sigma_q)(D)-q^{-8}(1-q^{-2})(\beta_{1/q}\circ \sigma_q)(A) & -q^{-6}( \alpha^-_{1/q}\circ \sigma_q)(B)
	\\[2pt]
	-q^{-6} (\alpha^-_{1/q}\circ \sigma_q) (C) & q^{-10}(\alpha^+_{1/q}\circ \sigma_q)(A)
\end{array}\!
\right).
\ee

Now we are ready to write down the characteristic identities (\ref{CHSp-1}) and (\ref{preCHSp}) for the case of
the standard $Sp(4)$-type RTT-algebra:
\ba
\lb{CH-RTT-Sp4}
M^{\overline{4}} - q\, M^{\overline{3}}\,\epsilon_1 +q^2 M^{\overline{2}}\,\epsilon_2- q^3 M \, \epsilon_1 g +q^4 I\, g^2 & =& 0,
\\[2pt]
\lb{preCH-RTT-Sp4}
M^{\overline{2}} - q\, M\, \epsilon_1 +q^2 I\, \epsilon_2 -q^3 \pi(M)\, \epsilon_1 +q^4 \pi^{\overline{2}}(M) & =& 0.
\ea
For reader's convenience we recall formulas for the powers of quantum matrices:
$$
M^{\overline{i+1} }= M\cdot \phi(M^{\overline{i}})\;\; \forall\, i\geq 1,\quad  \pi^{\overline{2}}(M)=-q^{-5}\, \phi^{-1}(\xi(M)\cdot \pi(M)),
$$
where in the last formula we took into account that $\mu =-q^{-5}$ and $G=I$ in our particular case.

\vskip .1cm
Using the definitions of the maps $\xi$, $\phi^{\pm 1}$, $\pi$ and of the elements 
$\epsilon_1$, $\epsilon_2$  given above, and applying the quadratic relations from the Appendix one can  check the parent characteristic identity (\ref{preCH-RTT-Sp4}) directly. The Cayley-Hamilton identity (\ref{CH-RTT-Sp4}) follows from it by
the $\star$-multiplication by $M^{\overline{2}}$.
\medskip

Finally, we consider the classical limit of the parent Cayley-Hamilton identities. In the limit $q\rightarrow 1$ the standard $Sp(2k)$-type R-matrix (\ref{R-Sp}) becomes the usual permutation and the 
quadratic relations (\ref{qma}) in the corresponding algebra ${\cal M}(P,P)$ imply the commutativity of the components of matrix $M$. The rank~$=1$ projector $K^{\mbox{\tiny (st)}}$ (\ref{K-Sp}) decouples from the R-matrix and the $g$-invariance conditions (\ref{tau2}) become independent of (\ref{qma}) and should be treated separately. We rewrite them in the familiar form
\be
\lb{g-inv-class}
M^t \,\Omega\, M\, =\, g\, \Omega \, =\, M\, \Omega\, M^t.
\ee
Here $M^t$ is the transposed matrix and $\Omega$ is the $2k\times 2k$ matrix of the symplectic quadratic form. With our choice of
the rank~$=1$ matrix $K^{\mbox{\tiny (st)}}$ it reads
$$
\Omega\, =\, 
\left(\!\!
\begin{array}{cc}
0 &  w
\\
-w & 0
\end{array}\!\!
\right),
$$
where $w$ is the $k\times k$ antidiagonal matrix: $w^i_j=\delta^i_{j'},\quad j'=k+1-j$. 

\vskip .1cm
Notice that in case $g\neq 0$ (more formally, if $g$ is invertible) the left and right equalities in (\ref{g-inv-class}) 
result in equivalent sets of conditions. 
On the contrary, in case $g=0$ these equalities are not equivalent and only together they give the complete set of 
the $g$-invariance conditions.

\vskip .1cm
Again, it is suitable to use the block notation (\ref{M-part}) for the matrix $M$, where now $A$, $B$, $C$ and $D$ 
are $k\times k$ matrices. The matrix $\pi(M)$ in this notation is
$$
\pi(M)\, =\, -\,\Omega\, M^t\, \Omega\, =\,
\left(\!\!
\begin{array}{cc}
D' &  -B'
\\
-C' & A'
\end{array}\!\!
\right), 
$$
where 
$X'= w X^t w$. This operation is a classical counterpart of the map $\sigma_q$ from our previous example.

\smallskip

The classical parent Cayley-Hamilton identity reads
\be
\lb{preCHSp-class}
\sum_{i=0}^{k} (-1)^i M^{k-i} \epsilon_i\, + \, \sum_{i=0}^{k-1} (-1)^{i} \pi(M)^{k-i} \epsilon_i\,=\,  0\, ,
\ee
where now all matrix powers are calculated according to the usual rules (the map $\phi$ in the classical limit is identical) and the coefficients $\epsilon_i$ become usual traces of the $i$-th wedge powers of the matrix $M$: $\epsilon_i = \tr (\wedge^i M)$ (the antisymmetrizers computed with the permutation matrix $P$ automatically include contributions from the 2-contraction $g$).
\smallskip

Assuming the invertibility of the matrix $A$ one can solve the $g$-invariance relations explicitly
\be
\lb{g-inv-solve}
M\, =\,
\left(\!\!
\begin{array}{cc}
	A &  A Y
	\\
	X A & X A Y +g A'^{-1}
\end{array}\!\!
\right)
\, =\, 
\left(\!\!
\begin{array}{cc}
	I &  0
	\\
	X & I
\end{array}\!\!
\right)
\left(\!\!
\begin{array}{cc}
	A & 0
	\\
	0 & g A'^{-1}
\end{array}\!\!
\right)
\left(\!\!
\begin{array}{cc}
	I &   Y
	\\
	0 & I
\end{array}\!\!
\right).
\ee
where matrices $X$, $Y$ are 
such that $X'=X, \;\; Y'=Y$. 

\vskip .1cm
Substituting this parameterization for $M$ into the identity  (\ref{preCHSp-class}) one can reduce it, 
 at  least in cases
$k=1,2,$ to the Cayley-Hamilton identities for  $k\times k$ matrices $A$ and $XY$. 

\subsection{Spectral parameterization}

\medskip
In this section we describe the parameterization
the coefficients of the characteristic polynomial
(\ref{CHSp-1}) by means of a $\Bbb C$-algebra
${\cal E}_{2k}$  of polynomials in $2k+1$ pairwise commuting variables
$\nu_i$, $~i=0,1,\dots ,2k$, satisfying conditions
\be
\lb{specSp}
\nu_{k+i}\, \nu_{k+1-i}\, =\, \nu_0^2\,\quad\forall\; i=1,2,\dots ,k\, .
\ee
We call $\nu_i$, $i=0,1,\dots ,2k$, {\em  spectral variables}. These variables play a role of the eigenvalues of the symplectic type quantum matrix $M$. This parameterization was initially aimed at comparing our results with
expressions given for the power sums for the RE-algebras in \cite{Mudr} (see
the subsection 8.3 there). Although the  derivation methods are very different the results agree up to some
obvious changes in a notation. Notice that compared to \cite{Mudr} we are working in a more general
setting. The generalization goes in several directions. First, we do not assume a ``standard''
Drinfel'd-Jimbo's form for the R-matrices defining the algebra and, moreover, we do not use any
deformation assumptions in our constructions. Next, we are working with a wider family of QM-algebras.
And, finally, we are working directly in the algebra without passing to representations\footnote{
Passing to the representations level is hardly possible except in the RE-algebra case. The reason is
that the characteristic subalgebra belongs to the center of the RE-algebra, which is not true for the
general QM-algebra.}.

\vskip .1cm
We are going to factorize the polynomial in the left hand side of the equation (\ref{CHSp-1}).
To this end, we realize elements of
the characteristic subalgebra ${\cal C}(R,F)$ as polynomials in the spectral variables and construct a
corresponding extention of the algebra ${\cal P}(R,F)$.

\begin{prop}
	\lb{corollary5.5}
In the setting of the theorem \ref{theorem5.4}, assume that the  elements $a_i$,
$i=1,2,\dots ,k$, are algebraically independent.
Consider an algebra homomorphism of the characteristic subalgebra ${\cal C}(R,F)$ to the algebra of
the spectral variables ${\cal E}_{2k}$,
$
\pi_{Sp(2k)}: {\cal C}(R,F)\rightarrow {\cal E}_{2k}\ ,
$
defined on the generators by
\be
\lb{rep-charSp}
\pi_{Sp(2k)}:\;\; g\mapsto \nu_0^2\, , \quad a_i\mapsto e_i(\nu_0,-\nu_0,\nu_1, \nu_2,\dots ,\nu_{2k})\,
\quad \forall\; i=1,\dots ,k\, ,
\ee
where $e_i$ are the elementary symmetric polynomials of their arguments (for the symmetric polynomials
we adopt a notation of \cite{Mac}). The map $\pi_{Sp(2k)}$ defines naturally a left
${\cal C}(R,F)$--module structure on the algebra ${\cal E}_{2k}$. Consider a corresponding completion
of the algebra ${\cal P}(R,F)$,
$$
{\cal P}_{Sp(2k)}(R,F)\, :=\, {\cal P}(R,F)\raisebox{-4pt}{$\bigotimes\atop
{\cal C}(R,F)$} {\cal E}_{2k}\, ,\vspace{-2mm}
$$
where the $\star \, $-product on the completed space is given by the formula
\be
\lb{P-Sp2}
(N\raisebox{-4pt}{$\bigotimes\atop {\cal C}(R,F)$}\nu)\star
(N'\raisebox{-4pt}{$\bigotimes\atop {\cal C}(R,F)$}\nu') :=
(N\star N')\raisebox{-4pt}{$\bigotimes\atop {\cal C}(R,F)$}(\nu\nu')\ \ \ \forall\;
N,N'\in {\cal P}(R,F)\ \; {\mathrm{and}}\ \; \forall\;\nu, \nu' \in {\cal E}_{2k}\, .
\ee
Then, in the completed algebra  ${\cal P}_{Sp(2k)}(R,F)$,
the Cayley-Hamilton identity (\ref{CHSp-1}) acquires a factorized form
\be
\lb{CHSp-factor} 
{\prod_{i=1}^{2k}}\hspace{-9.pt} {\scriptstyle \star}\; \left( M - q\nu_i I\right)\, =\, 0\, ,
\ee
where the  symbol $\displaystyle\prod\hspace{-10.4pt} {\scriptstyle \star} $ denotes the product with respect
to the $\star \, $-multiplication (\ref{P-Sp2}).
\end{prop}

\begin{rem}
{\rm For the classical symplectic groups, the functions $a_i$, $i=1,\dots ,k,$
on the manifold $Sp(2k)$ are
functionally independent. This justifies, at least perturbatively, the corresponding assumptions about the 
independence of the elements $a_i$ in the proposition above.}
\end{rem}

\begin{rem}
{\rm For a general  quantum matrix algebra ${\cal M}(R,F)$  the characteristic subalgebra does not belong to its 
center.  So, there is no general rule to define an extension by the spectral variables $\{\nu_i\}$ of the algebra 
${\cal M}(R,F)$. Nevertheless the commutative algebra ${\cal P}_{Sp(2k)}(R,F)$ admits the central extension by the spectral variables. Therefore we formulate the 
factorized Cayley-Hamilton identity for this extension.
	
\vskip .1cm	
However, for the reflection equation algebra ${\cal M}(R,R)$ the characteristic subalgebra lies in the center, the $\star$-product coincides with the usual matrix product and therefore one can assume that  eq.(\ref{CHSp-factor}) is satisfied in the central extension of
${\cal M}(R,R)$	by the spectral variables $\{\nu_i\}$.}
\end{rem}

\noindent {\bf Proof.}~ Using the equalities
$$
e_i(\nu_0,-\nu_0,\nu_1, \nu_2,\dots ,\nu_{2k})\, =\, e_i(\nu_1, \nu_2,\dots ,\nu_{2k}) -
\nu_0^2\, e_{i-2}(\nu_1, \nu_2,\dots ,\nu_{2k})\ \ \ \forall\ \ i\geq 0
$$
\vspace{-3mm}
and
$$
e_{k+i}(\nu_1,\nu_2,\dots ,\nu_{2k})\, =\,\nu_0^{2i}\, e_{k-i}(\nu_1,\nu_2,\dots ,
\nu_{2k})\ \ \ \forall\, i=1,\dots ,k, 
$$
if $\{\nu_i\}$ verifies eqs.(\ref{specSp}),
it is straightforward to check  that the map $\pi_{Sp(2k)}$ sends the coefficients (\ref{CHSp-2})
of the Cayley-Hamilton identity to the elementary symmetric functions in the spectral
variables: $\ \epsilon_i\ $ $\mapsto\ $
$e_i(\nu_1, \nu_2,\dots ,\nu_{2k})\ $ $\quad\forall\ i=1,\dots ,2k\, .$ \hfill$\blacksquare$
\medskip

In \cite{OP} we have derived the quantum analogs of the Newton and Wronsky relations among three series of 
elements of the characteristic subalgebra: the power sums $p_i$, the elementary symmetric functions $a_i$ and 
the complete symmetric functions $s_i$. Using these relations we now obtain the 
parameterization of the series $p_i$ and $s_i$ in terms of the spectral variables.  
\smallskip

\begin{prop}
\lb{corollary6.4}
Let  ${\cal M}(R,F)$ be the $Sp(2k)$-type  quantum matrix algebra. Assume that the algebra parameter $q$
fulfills the conditions $i_q\neq 0$, $i=2,\dots ,n$,  for some $n$.\footnote{For $n\leq k$ these conditions enter the initial settings for the $Sp(2k)$ type quantum matrix algebras.}
Then the elements $a_n$ and  $s_n$ can be defined recursively by the use of the Newton relations (see \cite{OP}, theorem 5.2)
\ba
\lb{Newton-a}
\sum_{i=0}^{n-1} (-q)^i a_i\, p_{n-i} &=& (-1)^{n-1} n_q\, a_n \, +\, (-1)^n
\sum_{i=1}^{\lfloor {n/2}\rfloor}\Bigl( \mu q^{n-2i} -q^{1-n+2i}\Bigr)\, a_{n-2i}\, g^i, 
\\[2pt]
\lb{Newton-s}
\sum_{i=0}^{n-1} q^{-i} s_i\, p_{n-i} &=&
  n_q\, s_n \, +\,\sum_{i=1}^{\lfloor {n/2}\rfloor}\Bigl(
\mu q^{2i-n} + q^{n-2i-1}\Bigr)\, s_{n-2i}\, g^i.
\ea
In this situation the elements $s_n$ and $p_n$ have the following images under the 
homomorphism $\pi_{Sp(2k)}$ (\ref{rep-charSp}):
\ba
\lb{para-s-Sp}
\pi_{Sp(2k)}:&& s_n \mapsto h_n(\nu_1,\nu_2,\dots ,\nu_{2k}) ,\qquad
p_n \mapsto q^{n-1}\sum_{i=1}^{2k} d_i \nu_i^n\, ,
\ea
where $h_n$ denotes the complete symmetric polynomial in its arguments and
\ba
\lb{para-p-Sp}
 d_i\, := {\nu_i - q^{-4}\nu_{2k+1-i}\over \nu_i -\nu_{2k+1-i}}
\prod _{j=1\atop j\neq i,\, 2k +1-i}^{2k} {\nu_i - q^{-2}\nu_j\over \nu_i - \nu_j}\, .
\ea

The power sums contain the rational functions $d_i$ in the spectral variables and are themselves
rational functions in $\{ \nu_i\}$. However, as it follows from the
Newton recursion (\ref{Newton-a}), the power sums simplify, modulo the
relations (\ref{specSp}), to polynomials in the spectral variables.\end{prop}

\noindent {\bf Proof.}~ For the proof, we use the following auxiliary statement:

\begin{lem}\lb{lemma6.5}
In the assumptions of proposition \ref{corollary6.4}, consider the iterations
\ba
\lb{mod-s}
s'_0=s_0\, , && s'_1=s_1\, , \quad\;\; s'_i = s_i + s'_{i-2}\, g\, ;
\\[2pt]
\lb{mod-p2} 
p'_0=(1-\mu^2 q^{2})/(q-q^{-1})\, ,  &&
p'_1=p_1\, , \quad\;\; p'_i = p_i + (q^{-2}p'_{i-2}-p_{i-2})\, g\;\;\; \forall i\geq 2.
\ea
The modified sequences 
$\{s'_i\}_{i=0}^n$,  $\{p'_i\}_{i=0}^n$ satisfy the following versions of the Newton and Wronski
relations
\ba
\lb{mod-N} 
\qquad \sum_{i=0}^{n-1} q^{-i} s_i p'_{n-i}& =&  n_q s_n\quad  \forall\; n\geq 1\, ;
\\
\lb{mod-W}
\sum_{i=0}^{n} (-1)^i a_i s'_{n-i}& =&\delta_{n,0}\,
\quad\; \forall\; n\geq 0\, .\ea\end{lem}

\noindent {\bf Proof.}~ For $n<2$, the equalities (\ref{mod-N})--(\ref{mod-W}) are clearly satisfied. For
$n\geq 2$, one can check them inductively, applying the iterative formulas (\ref{mod-s}), (\ref{mod-p2}).
\hfill$\blacksquare$

\medskip
We now notice that the images of the elements
$a_i$, $i=1,\dots ,n$, are given by the elementary symmetric functions (see eq.(\ref{rep-charSp})).
Hence, by the Wronski relations
(\ref{mod-W}), the images of the modified elements $s'_n$, $i=1,\dots ,n$, are the complete symmetric
functions in the same arguments. Using then  eq.(\ref{mod-s}) and taking into account the relation
$h_n(\nu_0,\nu_1,\dots)= \sum_{i=0}^n \nu_0^i\, h_{n-i}(\nu_1,\dots)$, it is easy to check the formulas
for the images of the elements $s_n$, which are given in eq.(\ref{para-s-Sp}).

\smallskip
To check the formulas for the power sums, we use the following statement, which was proved in \cite{GS}: if
the elements  $s_i$ for $i=0,1,\dots ,n\geq 1$ are realized as the complete
symmetric polynomials $h_i$
in some set of variables $\{\nu_i\}_{i=1}^{2k}$, then the elements $p'_n$,  defined
by eqs.(\ref{mod-N}), have the following expressions in terms of the variables $\nu_i$
\be
\lb{formula-GS}
 p'_n   = q^{n-1} \sum_{i-1}^p \widehat d_i \nu_i^n ,
\qquad \mbox{where}\quad \widehat d_i :=
\prod_{j=1\atop j\neq i}^{p}{\nu_i -q^{-2} \nu_j\over \nu_i - \nu_j}.
\ee

The proof of  
(\ref{para-s-Sp}), (\ref{para-p-Sp}) for the power sums $p_n$ proceeds  as follows.
\smallskip

Assuming that the relation (\ref{formula-GS}) stays valid for $p'_0$ (note, $p'_0$ is not fixed by the
recursion (\ref{mod-N})) and making the Ansatz (\ref{para-s-Sp}) for the power sums $p_i$ for
$i=0,1,\dots ,n$, we make use of the recursion (\ref{mod-p2}). Upon substitutions, we find that
the relations (\ref{mod-p2}) hold valid provided that
\be\lb{d-hat-d}d_i \, =\, {\nu_i^2 - q^{-4}\nu_0^2\over\nu_i^2 -q^{-2}\nu_0^2}\,\,\widehat d_i\, .\ee
Taking into account the relations (\ref{specSp}) for the spectral variables $\nu_i\in {\cal E}_{2k}$, we observe that the conditions (\ref{d-hat-d}) dictate the choice (\ref{para-p-Sp})
for $d_i$.

\smallskip
It remains to verify the initial settings for the recursion (\ref{mod-p2}). They are:
\ba
\lb{init-1}
p'_0& =&  q^{-1} \displaystyle \sum_{i=1}^{2k} \widehat d_i\,
=\, {1-\mu^2 q^2\over q-q^{-1}}|_{\mu=-q^{-1-2k}}\, =\, q^{-2k} (2k)_q\, ,
\\[2pt]
\lb{init-3}
p_1 &=&p'_1\quad \Leftrightarrow\quad \displaystyle
\sum_{i=1}^{2k}\nu_i (d_i - \widehat d_i)\, =\, 0, 
\ea
as well as the expression (\ref{P-i}) for $p_0$:
\ba
\lb{init-2}
p_0 & =&  q^{-1} \displaystyle \sum_{i=1}^{2k} d_i\,
=\, {\rm Tr}_R I|_{\mu=-q^{-1-2k}}\, =\, q^{-1-2k}\bigl( (2k+1)_q -1\bigr).
\ea

\smallskip
To verify them, we use expansions of the following rational functions
$$
w_1(z) := \prod_{i=1}^{2k} {z-q^{-2} \nu_i\over z-\nu_i}\, , \qquad
w_2(z) := {\nu_0^2 w_1(z)\over z^2 - q^{-2} \nu_0^2}\, , \qquad w_3(z) := z w_2(z)
$$
in simple ratios. 

Expanding $w_1(z)$ and evaluating the result at $z=0$, we prove immediately the
condition (\ref{init-1}).

A less trivial check of
the condition (\ref{init-2}) we comment in more details. Expanding $w_2(z)$, we obtain
$$
w_2(z)\, =\, \sum_{i=1}^{2k} q^2 (d_i-\widehat d_i) {\nu_i\over z-\nu_i}\, +\,  {q\nu_0\over 2}\Bigl(
{w_1(q^{-1}\nu_0)\over z-q^{-1}\nu_0} -{w_1(-q^{-1}\nu_0)\over z+q^{-1}\nu_0}\Bigr)\, .
$$
Here, for the transformation of the first term in the right hand side, we used the formulas
(\ref{formula-GS}) and (\ref{d-hat-d}) and applied the relations (\ref{specSp}), which confine the variables
$\nu_i\in {\cal E}_{2k}$. The  relations (\ref{specSp}) also allow us to calculate
$w_1(\pm q^{-1}\nu_0)=q^{-2k}$. Thus, evaluating $w_2(z)$ at $z=0$, we obtain
$$
w_2(0)\, =\, -q^{2-4k}\, =\,  -q^3 (p_0-p'_0)\, -\, q^{2-2k}\, ,
$$
wherefrom the condition (\ref{init-2}) follows.

\vskip .1cm
A check of the condition (\ref{init-3}), by the expansion
and evaluation of $w_3(z)$ at $z=0$, is a similar calculation. \hfill$\blacksquare$

\appendix
\section{Standard  $Sp(4)$-type RTT-algebra.}
\lb{Sp(4)}

\def\theequation{\thesection.\arabic{equation}}
\makeatletter\@addtoreset{equation}{section}\makeatother

Here we present quadratic relations for the $Sp(4)$-type RTT-algebra ${\cal M}(R,P)$ corresponding to the standard symplectic R-matrix
(\ref{R-Sp}), where we take $k=2$. The 2-contraction $g$ turns out to be central in this algebra.
With the additional condition $g=\mu^2\, 1=q^{-10}\, 1$ this algebra can be interpreted as the quantized algebra of functions on the Lie group $Sp(4)$.
\smallskip

For the $4\times 4$ matrix of generators of this algebra we use the following notation
\be
\lb{Sp4-M}
M\, =\,
\left(
\begin{array}{cc}
A & B
\\
C & D
\end{array}
\right),
\ee
where $A$, $B$, $C$, and $D$ are $2\times 2$ matrices.  For any matrix $X\in\{ A,B,C,D\}$ we denote its matrix 
components as $X^{i}_{j}$, $i,j\in\{1,2\}$.

\vskip .1cm
Quadratic relations in the standard $Sp(4)$-type RTT-algebra contain 120 permutation relations for 16 
components of the quantum matrix $M$ and
10 additional conditions (\ref{tau2}) which are responsible for  
invariance of the symplectic form encoded in the 
rank$=1$ projector $K$.
\smallskip

For the presentation of the permutation relations we fix the following linear order on the  components of $M$:
$$
X^{1}_{1}<X^{1}_{2}<X^{2}_{1}<X^{2}_{2}\;\; \forall \, X\in\{A,B,C,D\}, 
\quad D^{i_1}_{j_1}< C^{i_2}_{j_2} < B^{i_3}_{j_3}< A^{i_4}_{j_4}\;\; \forall\, i_k,j_k=\{1,2\}.
$$
Permutation relations among the components of the 2x2 matrices $A$, \dots ,$D$ take a universal form:
\ba
\nn
q\, X^{i}_{2} X^{i}_{1}\, =\, X^{i}_{1} X^{i}_{2}, \quad  q\, X^{2}_{i} X^{1}_{i}\, =\, X^{1}_{i} X^{2}_{i},
\quad [X^{2}_{1}, X^{1}_{2}] \, =\, 0,\quad
[X^{2}_{2}, X^{1}_{1}] \, =\, -\lambda\, X^{1}_{2} X^{2}_{1},
\ea
where in the last formula and below we use the shorthand notation $\lambda := q-q^{-1}$.

The rest 96 permutation relations between the components of different matrices $A$, $B$, $C$ and $D$ are separated into eight  groups according to the type of permutation. In the formulas below 
indices $i$, $j$ 
take values 1 or 2;~ $i':=3-i$.

\vskip .1cm
Commutators:
\ba
\nn
[A^{2}_{i},B^{1}_{i}]&=& [A^{i}_{2},C^{i}_{1}]\, =\,  [B^{i}_{2}, D^{i}_{1}]\, =\, [C^{2}_{i}, D^{1}_{i}]\, =\, 0,
\\[2pt]
\nn
[B^{i}_{j}, C^{i}_{j}]& =& 0, \qquad\qquad\;\; [B^{1}_{2}, C^{2}_{1}]\, =\, 0;
\ea
$q$-commutators:
\ba
\nn
A^{i}_{j}B^{i}_{j}- q\, B^{i}_{j}A^{i}_{j}&=& A^{i}_{j}C^{i}_{j}- q\, C^{i}_{j}A^{i}_{j}\, =\, B^{i}_{j}D^{i}_{j}- q\, D^{i}_{j}B^{i}_{j}\, =\, C^{i}_{j}D^{i}_{j}- q\, D^{i}_{j}C^{i}_{j}\, =\, 0,
\\[2pt]
\nn
A^{2}_{1}B^{1}_{2}- q\, B^{1}_{2}A^{2}_{1}&=& A^{1}_{2}C^{2}_{1}- q\, C^{2}_{1}A^{1}_{2}\, =\, B^{1}_{2}D^{2}_{1}- q\, D^{2}_{1}B^{1}_{2}\, =\, C^{2}_{1}D^{1}_{2}- q\, D^{1}_{2}C^{2}_{1}\, =\, 0,
\\[2pt]
\nn
B^{1}_{i}C^{2}_{i}- q\, C^{2}_{i}B^{1}_{i}&=&  0, \qquad\qquad\;\; B^{i}_{2}C^{i}_{1}- q^{-1} C^{i}_{1}B^{i}_{2}\, =\,  0;
\ea
$q^2$-commutators:
\ba
\nn
A^{i}_{1}B^{i}_{2}- q^2 B^{i}_{2}A^{i}_{1}&=& A^{1}_{i}C^{2}_{i}- q^2 C^{2}_{i}A^{1}_{i}\, =\, B^{1}_{i}D^{2}_{i}- q^2 D^{2}_{i}B^{1}_{i}\, =\, C^{i}_{1}D^{i}_{2}- q^2 D^{i}_{2}C^{i}_{1}\, =\, 0;
\ea
commutators with $\pm\lambda$-additional term (this just means that the numeric coefficient of an extra term
is equal to $\pm\lambda$):
\ba
\nn
[A^{1}_{i},B^{2}_{i}]-\lambda\, B^{1}_{i}A^{2}_{i}\!\!\!&=&\!\!\! [A^{i}_{1},C^{i}_{2}]-\lambda\, C^{i}_{1}A^{i}_{2}\, =\,  [B^{i}_{1},D^{i}_{2}]-\lambda\, D^{i}_{1}B^{i}_{2}\, =\,
[C^{1}_{i},D^{2}_{i}]-\lambda\, D^{1}_{i}C^{2}_{i}\, =\, 0,
\\[2pt]
\nn
[A^{i}_{j}, D^{i}_{j}]-\lambda\, C^{i}_{j} B^{i}_{j}\!\!\!&=&\!\!\! 0, \qquad\quad [B^{2}_{2}, C^{1}_{1}]-\lambda\, C^{2}_{1} B^{1}_{2}\, =\, 0,
\qquad\quad [B^{1}_{1}, C^{2}_{2}]+\lambda\, C^{2}_{1} B^{1}_{2}\, =\, 0;
\ea
$q$-commutators with $\pm q^{\pm 1}\lambda$-additional term: 
\ba
\nn
A^{i}_{i}B^{i'}_{i'}- q\, B^{i'}_{i'}A^{i}_{i}- \lambda q\, B^{1}_{2}A^{2}_{1}& =& A^{i}_{i}C^{i'}_{i'}- q\, C^{i'}_{i'}A^{i}_{i}- \lambda q\, C^{2}_{1}A^{1}_{2}\, =\, 0,
\\[2pt]
\nn
B^{i}_{i}D^{i'}_{i'}- q\, D^{i'}_{i'}B^{i}_{i}- \lambda q\, D^{2}_{1}B^{1}_{2}& =& C^{i}_{i}D^{i'}_{i'}- q\, D^{i'}_{i'}C^{i}_{i}- \lambda q\, D^{1}_{2}C^{2}_{1}\, =\, 0,
\\[2pt]
\nn
A^{1}_{i}D^{2}_{i}- q\, D^{2}_{i}A^{1}_{i}- \lambda q\, C^{2}_{i}B^{1}_{i}&=& B^{2}_{i}C^{1}_{i}- q\, C^{1}_{i}B^{2}_{i}- \lambda q\, C^{2}_{i}B^{1}_{i}\, =\, 0,
\ea
\vspace{-6.5mm}
\ba
\nn
\phantom{a}\mbox{\hspace{8mm}}B^{i}_{1}C^{i}_{2}- q^{-1} C^{i}_{2}B^{i}_{1}+ \lambda q^{-1} C^{i}_{1}B^{i}_{2}& =& 0;
\ea
$q$-commutators with $\pm \lambda$-additional term: 
\ba
\nn
A^{i}_{1}D^{i}_{2}- q\, D^{i}_{2}A^{i}_{1}- \lambda\,  C^{i}_{1}B^{i}_{2}& =&  0.
\ea
$q^2$-commutators with $\pm q^2\lambda$-additional term: 
\ba
\nn
A^{i}_{2}B^{i}_{1}- q^2 B^{i}_{1}A^{i}_{2}- \lambda q^2 B^{i}_{2}A^{i}_{1}& =& A^{2}_{i}C^{1}_{i}- q^2 C^{1}_{i}A^{2}_{i}- \lambda q^2 C^{2}_{i}A^{1}_{i}\, =\, 0,
\\[2pt]
\nn
B^{2}_{i}D^{1}_{i}- q^2 D^{1}_{i}B^{2}_{i}- \lambda q^2 D^{2}_{i}B^{1}_{i}& =& C^{i}_{2}D^{i}_{1}- q^2 D^{i}_{1}C^{i}_{2}- \lambda q^2 D^{i}_{2}C^{i}_{1}\, =\, 0.
\ea
More complicated relations:
\ba
\nn
A^{i}_{2}D^{i}_{1}- q^{-1} D^{i}_{1}A^{i}_{2}&=& \lambda q^{-3} C^{i}_{1}B^{i}_{2}+\lambda  2_q q^{-1} C^{i}_{2}B^{i}_{1},
\\[2pt]
\nn
A^{2}_{i}D^{1}_{i}- q^{-1} D^{1}_{i}A^{2}_{i}&=& \lambda q^{2} C^{2}_{i}B^{1}_{i}+\lambda 2_q\, C^{1}_{i}B^{2}_{i},
\\[2pt]
\nn
A^{1}_{2}B^{2}_{1}- q^{-1} B^{2}_{1}A^{1}_{2}&=& \lambda q^{2} B^{1}_{2}A^{2}_{1}+\lambda 2_q\, B^{1}_{1}A^{2}_{2},
\\[2pt]
\nn
A^{2}_{1}C^{1}_{2}- q^{-1} C^{1}_{2}A^{2}_{1}&=& \lambda q^{2} C^{2}_{1}A^{1}_{2}+\lambda 2_q\, C^{1}_{1}A^{2}_{2},
\\[2pt]
\nn
B^{2}_{1}D^{1}_{2}- q^{-1} D^{1}_{2}B^{2}_{1}&=& \lambda q^{2} D^{2}_{1}B^{1}_{2}+\lambda 2_q\, D^{1}_{1}B^{2}_{2},
\\[2pt]
\nn
C^{1}_{2}D^{2}_{1}- q^{-1} D^{2}_{1}C^{1}_{2}&=& \lambda q^{2} D^{1}_{2}C^{2}_{1}+\lambda 2_q\, D^{1}_{1}C^{2}_{2},
\ea
\vspace{-9mm}
\ba
\nn
[B^{2}_{1},C^{1}_{2}]&=& \lambda\,  C^{2}_{2}B^{1}_{1}-\lambda\, C^{1}_{1}B^{2}_{2}+\lambda^2 C^{2}_{1} B^{1}_{2},
\\[2pt]
\nn
[A^{1}_{1},D^{2}_{2}]&=& -\lambda\,  D^{1}_{2}A^{2}_{1}+\lambda q^{-2} C^{1}_{1}B^{2}_{2}+\lambda 2_q\, C^{2}_{1} B^{1}_{2},
\\[2pt]
\nn
[A^{1}_{2},D^{2}_{1}]&=& \lambda q^{-2} C^{2}_{1}B^{1}_{2}+\lambda 2_q\, C^{2}_{2} B^{1}_{1},
\\[2pt]
\nn
[A^{2}_{1},D^{1}_{2}]&=& \lambda q^{2} C^{2}_{1}B^{1}_{2}+\lambda 2_q\, C^{1}_{1} B^{2}_{2},
\ea
\vspace{-9mm}
\ba
\nn
[A^{2}_{2},D^{1}_{1}]\!\!\!&=&\!\!\! \lambda\,  D^{2}_{1}A^{1}_{2}+\lambda q^{-2} C^{1}_{1}B^{2}_{2}+\lambda 2_q\, C^{1}_{2} B^{2}_{1}+\lambda^2 q^{-2} C^{2}_{1}B^{1}_{2}+
\lambda^2 2_q\, C^{2}_{2}B^{1}_{1}.
\ea

Finally, from the matrix relations
$$
M_1 M_2 K_{12} = K_{12} M_1 M_2 = \mu^{-2} K_{12}\, g = q^{10} K_{12}\, g
$$
we extract two equivalent expressions for $g$
\ba
\nn
g & =& q^{-10}\left( D^{1}_{1} A^{2}_{2} + q D^{1}_{2} A^{2}_{1}-q^{-2} C^{1}_{2} B^{2}_{1} - q^{-3} C^{1}_{1} B^{2}_{2} \right)
\\
\lb{g-RTT-Sp4}
&=& q^{-10}\left(D^{2}_{2} A^{1}_{1} +q^{-1} D^{1}_{2} A^{2}_{1} - q^{-2} C^{2}_{1} B^{1}_{2} -q^{-3} C^{1}_{1} B^{2}_{2}\right),
\ea
and 10 invariance conditions 
\ba
\nn
B^{1}_{1} A^{2}_{2} + q\, B^{1}_{2} A^{2}_{1}-q\, B^{2}_{1} A^{1}_{2} - q^2 B^{2}_{2} A^{1}_{1}& =&
D^{1}_{1} C^{2}_{2} + q\, D^{1}_{2} C^{2}_{1}-q\, D^{2}_{1} C^{1}_{2} - q^2 D^{2}_{2} C^{1}_{1}
\,=\,0,
\\[2pt]
\nn
C^{1}_{1} A^{2}_{2} + q\, C^{2}_{1} A^{1}_{2}-q\, C^{1}_{2} A^{2}_{1} - q^2 C^{2}_{2} A^{1}_{1}& =&
D^{1}_{1} B^{2}_{2} + q\, D^{2}_{1} B^{1}_{2}-q\, D^{1}_{2} B^{2}_{1} - q^2 D^{2}_{2} B^{1}_{1}
\,=\,0,
\\[2pt]
\nn
C^{i}_{1} B^{i}_{2} + q\, C^{i}_{2} B^{i}_{1}-q^3 D^{i}_{1} A^{i}_{2} - q^4 D^{i}_{2} A^{i}_{1}& =&
C^{1}_{i} B^{2}_{i} + q\, C^{2}_{i} B^{1}_{i}-q\, D^{1}_{i} A^{2}_{i} - q^2 D^{2}_{i} A^{1}_{i}
\,=\,0,
\ea
\vspace{-8mm}
\ba
\nn
C^{1}_{1} B^{2}_{2} - C^{2}_{2} B^{1}_{1} +\lambda\, C^{2}_{1} B^{1}_{2} -q^2 D^{1}_{2} A^{2}_{1} +q^2 D^{2}_{1} A^{1}_{2} &=&0,
\\[2pt]
\nn
C^{1}_{2} B^{2}_{1} - C^{2}_{1} B^{1}_{2}  -q^2 D^{1}_{1} A^{2}_{2} +q^2 D^{2}_{2} A^{1}_{1} -\lambda q^2  D^{1}_{2} A^{2}_{1}&=&0.
\ea

\section*{Acknowledgments}

The work of the first author (O. O.) was supported by the grant PhyMath ANR-19-CE40-0021, 
the Program of Competitive Growth of Kazan Federal University and the Russian-French Poncelet Laboratory
UMI 2615.
The work of the second author (P. P.) was  supported
by the Academic Fund Program at the HSE University (grant no.20-01-032 for
the years 2020-2022) and the Russian Academic Excellence Project `5-100',
and by the grant of RFBR no.20-51-12005.

\bigskip
\addtocontents{toc}{\contentsline {section}{\numberline {} References}{\pageref{refer}}}

\end{document}